\documentclass[10pt,a4paper]{amsart}
\usepackage{verbatim}

\usepackage[toc]{appendix}
\usepackage[T1]{fontenc}
\usepackage{graphicx}
\usepackage{enumerate}
\usepackage{amsmath,amsfonts,amssymb}
\usepackage{color}
\def\loc{\operatorname{loc}}
\usepackage{cite}
\definecolor{citation}{rgb}{0.2,0.58,0.2} 
\definecolor{formula}{rgb}{0.1,0.2,0.6}
\definecolor{url}{rgb}{0.3,0,0.5}
\usepackage{pgf,tikz}
\usepackage{mathrsfs}
\usepackage{fancyhdr}
\usepackage{dutchcal}

\usepackage{dsfont}

\newcommand{\reqnomode}{\tagsleft@false}

\usepackage[colorlinks,pdfpagelabels,pdfstartview = FitH,bookmarksopen = true,bookmarksnumbered = true,urlcolor=url,linkcolor = formula,plainpages = false,hypertexnames = false,citecolor = citation] {hyperref}

\def\dx{\,{\rm d}x}

\def\dtt{\,{\rm d}t}
\def\dll{\,{\rm d}l}
\def\ds{\,{\rm d}s}
\def\dt{\,{\rm d}t}
\def\dy{\,{\rm d}y}
\def\dzz{\,{\rm d}z'}
\def \d{\,{\rm d}}
\def \diver{\,{\rm div}}
\def\dist{\,{\rm dist}}

\def\supp{\,{\rm supp}}

\allowdisplaybreaks
\makeatletter
\DeclareRobustCommand*{\bfseries}{%
  \not@math@alphabet\bfseries\mathbf
  \fontseries\bfdefault\selectfont
  \boldmath
}

\makeatother

\newlength{\defbaselineskip}
\newcommand{\setlinespacing}[1]
           {\setlength{\baselineskip}{#1 \defbaselineskip}}
\newcommand{\mint}{\mathop{\int\hskip -1,05em -\, \!\!\!}\nolimits}

\newtheorem{theorem}{Theorem}

\newtheorem{definition}{Definition}
\newtheorem{remark}{Remark}[section]
\newtheorem{lemma}{Lemma}[section]
\newtheorem{proposition}{Proposition}[section]
\numberwithin{equation}{section}
\allowdisplaybreaks[4]
\newcommand{\ta}{a_{j}}

\newcommand{\kk}{\kappa}

\def\er{\mathbb R}

\newcommand{\ti}[1]{\tilde{#1}}

\newcommand\eps\varepsilon
\def\eqn#1$$#2$${\begin{equation}\label#1#2\end{equation}}

\newcommand{\be}{\begin{equation}}
\newcommand{\ee}{\end{equation}}

\newcommand{\rr}{\varrho}

\newcommand{\snr}[1]{\lvert #1\rvert}

\newcommand{\nr}[1]{\lVert #1 \rVert}

\newcommand{\N}{\mathbb{N}}

\def\name[#1, #2]{#1 #2}

\newcommand{\Di}{\Delta_{h}}
\newcommand{\taa}{\tilde{a}}

\title[Parabolic equations with $(p,q)$-growth]{Gradient bounds for solutions to irregular parabolic equations with $(p,q)$-growth}

\author[De Filippis]{Cristiana De Filippis}  \address{Cristiana De Filippis\\Mathematical Institute, University of Oxford\\ Andrew Wiles Building, Radcliffe Observatory Quarter, Woodstock Road, Oxford, OX26GG, Oxford, United Kingdom} \email{\texttt{Cristiana.DeFilippis@maths.ox.ac.uk}}
\begin{document}

\subjclass[2010]{35K20, 35K61, 35K92\vspace{1mm}} %%ALERT CHECK 35J60 23J70 35B65 35D40

\keywords{Regularity, parabolic equations, $(p,q)$-growth\vspace{1mm}}

\thanks{{\it Acknowledgements.}\ This work is supported by the Engineering and Physical Sciences Research Council (EPSRC): CDT Grant Ref. EP/L015811/1. 
\vspace{1mm}}

\maketitle

\begin{abstract}
We provide quantitative gradient bounds for solutions to certain parabolic equations with unbalanced polynomial growth and non-smooth coefficients.
\end{abstract}
\vspace{3mm}
{\small \tableofcontents}

\setlinespacing{1.08}
\section{Introduction}
We focus on the Cauchy-Dirichlet problem
\begin{flalign}\label{pdd}
\begin{cases}
\ \partial_{t}u-\diver \ a(x,t,Du)=0\quad &\mbox{in} \ \ \Omega_{T}\\
\ u=f\quad &\mbox{on} \ \ \partial_{par}\Omega_{T},
\end{cases}
\end{flalign}
with initial-boundary datum $f\colon \mathbb{R}^{n+1}\to \mathbb{R}$ as in \eqref{ggg} below and nonlinear diffusive tensor $a(\cdot)$ featuring $(p,q)$-growth conditions as displayed in \eqref{ref}. The main novelties here are twofold: the map $x\mapsto a(x,t,z)$ is only Sobolev-differentiable in the sense that 
\begin{flalign*}
\snr{\partial_{x}a(x,t,z)}\le \gamma(x,t)\left[(\mu^{2}+\snr{z}^{2})^{\frac{p-1}{2}}+(\mu^{2}+\snr{z}^{2})^{\frac{q-1}{2}}\right],
\end{flalign*}
where $\gamma$ possess a suitably high degree of integrability, cf. \eqref{gamma}. Moreover, we can treat in a single shot both the degenerate case $p\ge 2$ and the singular one $p<2$, allowing also for the case $\mu=0$. Precisely, we prove that
\begin{theorem}\label{t1}
If assumptions \eqref{regg}-\eqref{ggg} are satisfied, Cauchy-Dirichlet problem \eqref{pdd} admits a solution $u\in L^{p}(0,T;W^{1,p}(\Omega))$ such that
\begin{flalign}\label{t1.1}
Du\in L^{\infty}_{\loc}(\Omega_{T},\mathbb{R}^{n}),\qquad V_{\mu,p}(Du)\in L^{2}_{\loc}(0,T;W^{1,2}_{\loc}(\Omega,\mathbb{R}^{n}))
\end{flalign}
and 
\begin{flalign}\label{t1.4}
u\in W^{\iota,2}_{\loc}(0,T;L^{2}_{\loc}(\Omega))\quad \mbox{for all} \ \ \iota\in \left(0,\frac{1}{2}\right).
\end{flalign}
In particular, if $Q_{\rr}\Subset \Omega_{T}$ is any parabolic cylinder there holds that
\begin{flalign}\label{t1.5}
\nr{H(Du)}_{L^{\infty}(Q_{\rr/2})}\le \frac{c}{\rr^{\beta_{1}}}\left[1+\left(\mint_{Q_{\rr}}H(Du)^{\frac{p}{2}} \ \dy\right)^{\beta_{2}}\right],
\end{flalign}
with $c\equiv c(\texttt{data})$ and $\beta_{1},\beta_{2}\equiv \beta_{1},\beta_{2}(n,p,q,d)$.
\end{theorem}
We refer to Sections \ref{not}-\ref{ma} for a detailed description of the various quantities involved in the previous statement. Our analysis includes equations with double phase structure, such as
\begin{flalign*}
&\partial_{t}u-\diver\left(\snr{Du}^{p-2}Du+b(x,t)\snr{Du}^{q-2}Du\right)=0\quad \mbox{in} \ \ \Omega_{T}\\
&b\in L^{\infty}(\Omega_{T})\ \  \mbox{with}\ \ \partial_{x}b\in L^{d}(\Omega_{T}); 
\end{flalign*}
or equations with variable exponent:
\begin{flalign*}
&\partial_{t}u-\diver\left(\snr{Du}^{p(x,t)-2}Du\right)=0\quad \mbox{in} \ \ \Omega_{T}\\
&p\in L^{\infty}(\Omega_{T})\ \  \mbox{with}\ \ \partial_{x}p\in L^{d}(\Omega_{T});
\end{flalign*}
and also anisotropic equations like
\begin{flalign*}
&\partial_{t}u-\left[\diver(\snr{Du}^{p-2}Du)+\sum_{i=1}^{n}\partial_{x_{i}}\left((\mu^{2}+\snr{\partial_{x_{i}}u}^{2})^{\frac{p_{i}-2}{2}}\partial_{x_{i}}u\right)\right]\\
&1<p\le p_{i}<\infty \ \ \mbox{for all} \ \ i\in \{1,\cdots,n\},
\end{flalign*}
where $(p,q)$, $(\inf_{\Omega_{T}}p,\sup_{\Omega_{T}}p)$, $\left(p,\max_{i\in \{1,\cdots,n\}}p_{i}\right)$ satisfy \eqref{pq} and $d$ is described by $\eqref{gamma}_{2}$. To the best of our knowledge, the result stated in Theorem \ref{t1} is new already in the standard growth case $p=q$. This fact poses additional difficulties due to the lack of informations on the regularity of solutions to \eqref{pdd} when $a(\cdot)$ has balanced polynomial growth. To overcome this issue, we proceed in two steps: first, we prove an higher integrability result for solutions to a regularized version of problem \eqref{pdd}. Then we use it to construct a sequence of maps satisfying suitable uniform estimates and converging to a solution of \eqref{pdd}. For the sake of simplicity, Theorem \ref{t1} is proved in the scalar case, but all our arguments can be adapted in a straightforward way to the vectorial setting as well, provided that $a(\cdot)$ has radial structure. Let us put our result in the context of the available literature. The systematic study of problem
\begin{flalign}\label{elc}
\begin{cases}
\ -\diver \ a(x,Du)=0\quad &\mbox{in} \ \ \Omega\\
\ u=f\quad &\mbox{on} \ \ \partial\Omega
\end{cases}
\end{flalign}
i.e., the elliptic counterpart of \eqref{pdd} started in \cite{ma2,ma3,ma1} and, subsequently, has undergone an intensive development over the last years, see for instance \cite{bacomi,bemi,besc,besc1,bobr,dm1,dm,deoh,dima,sharp,haok,lieb} and references therein. As suggested by the counterexamples contained in \cite{sharp,ma2}, already in the elliptic setting the regularity of solution to \eqref{elc} is strongly linked to the closeness of the exponents $(p,q)$ ruling the growth of the vector field $a(\cdot)$. Precisely, it turns out that
\begin{flalign}\label{pq.2}
1\le\frac{q}{p}<1+\mathcal{M}(\texttt{problem's data}),
\end{flalign}
where $\mathcal{M}(\cdot)$ is in general a bounded function connecting the various informations given \emph{a priori} about solutions. In this respect, we refer to \cite{bacomi} for an idea on the subtle yet quantifiable interplay between the regularity of solutions and the main parameters of the problem and to \cite{bemi,dm1,dm}, where is shown that, as long as $p$ and $q$ stay close to each other, problems with $(p,q)$-growth can be interpreted as perturbations of problems having standard $p$-growth. In the parabolic setting, the regularity for solutions of \eqref{pdd} is very well understood when $a(\cdot)$ is modelled upon the parabolic $p$-laplacian, see e.g. \cite{d,dumi1,dumi,km2,km3} for an overview of the state of the art on this matter and \cite{ba,ba1}, where more general structures are analyzed. Finally, the question of existence of regular solutions of \eqref{pdd} when the nonlinear tensor $a(\cdot)$ has unbalanced polynomial growth was treated in \cite{bdm,bdm1,si,s}. The theory exposed in these papers confirms that, as in the elliptic case, a restriction like \eqref{pq.2} on the ratio $q/p$ suffices to prove existence of regular solutions to \eqref{pdd}. Actually, the function $\mathcal{M}(\cdot)$ is worsen for parabolic equations than for elliptic ones, due to the so-called phenomenon of caloric deficit, originated from the difference of scaling in space and time, see e.g. \cite{bdm,s}, in which $\mathcal{M}(\cdot)$ is quantified as a function of $n$ and $p$. In our case, $\mathcal{M}(\cdot)$ has to take into account also the integrability exponent of $\gamma$, therefore it depends on $n,p,d$ and, reversing the process of caloric deficit, it renders precisely the bound for elliptic equations with Sobolev-differentiable coefficients appearing in \cite{dm1,dm,elemarmas}.
\subsubsection*{Organization of the paper} This paper is organized as follows: in Section \ref{pre} contains our notation, the list of the assumptions which will rule problem \eqref{pdd} and several by now classical tools in the framework of regularity theory for elliptic and parabolic PDE. Sections \ref{high}-\ref{mose} are devoted to the proof of Proposition \ref{p1} and Theorem \ref{t1} respectively.  

\section{Preliminaries}\label{pre}
In this section we display the notation adopted throughout the paper and list some well-known result which will be helpful in the various proofs presented.
\subsection{Notation}\label{not}
In this paper, $\Omega_{T}:=\Omega\times (0,T)$ is a space-time cylinder over an open, bounded domain $\Omega\subset \mathbb{R}^{n}$, $n\ge 2$ with $C^{1}$-boundary. If $\tilde{\Omega}\subseteq \Omega$ and $t_{0}\in [0,T]$, by $\Omega_{t_{0}}$ we mean the subcylinder $\Omega\times (0,t_{0})\subseteq \Omega_{T}$. Clearly, when $t_{0}=0$, $\Omega_{0}\equiv \Omega$. We denote by $B_{\rr}(x_{0}):=\left\{x\in \mathbb{R}^{n}\colon \snr{x-x_{0}}<\rr\right\}$ the $n$-dimensional open ball centered at $x_{0}\in \mathbb{R}^{n}$ with radius $\rr>0$. When working in the parabolic setting it is convenient to consider parabolic cylinders
\begin{flalign*}
Q_{\rr}(y_{0}):=B_{\rr}(x_{0})\times (t_{0}-\rr^{2},t_{0})\quad \mbox{where} \ \ y_{0}:=(x_{0},t_{0})\in \mathbb{R}^{n+1},
\end{flalign*}
i.e., balls in the parabolic metric. With "$y$" we shall always denote the couple $(x,t)\in \Omega_{T}$. Very often, when not otherwise stated, different balls (or cylinders) in the same context will share the same center. Given any differentiable map $G: \Omega\times \mathbb{R}\times \mathbb{R}^{n}\to \mathbb{R}$, with $\partial_{z}G(x,t,z)$ we mean the derivative of $G(\cdot)$ with respect to the $z$ variable, by $\partial_{t}G(x,t,z)$ the derivative in the time variable $t$ and by $\partial_{x}G(x,t,z)$ the derivative of $G$ with respect to the space variable $x$. We name "$c$" a general constant larger than one. Different occurrences from line to line will be still denoted by $c$, while special occurrences will be denoted by $c_1, c_2,  \tilde c$ and so on. Relevant
dependencies on parameters will be emphasized using parentheses, i.e., $c_{1}\equiv c_1(n,p)$ means that $c_1$ depends on $n,p$. For the sake of clarity, we shall adopt the shorthand notation
\begin{flalign*}
\texttt{data}:=\left(n,\nu,L,p,q,d,\nr{\gamma}_{L^{d}(\Omega_{T})}\right).
\end{flalign*}
In most of the inequalities appearing in the proof of our results we will use the symbols "$\lesssim$" or "$\gtrsim$", meaning that the inequalities hold up to constants depending from some (or all) the parameters collected in $\texttt{data}$. We refer to Section \ref{ma} for more details on the quantities appearing in the expansion of $\texttt{data}$.
\subsection{Main assumptions}\label{ma}
When dealing with the Cauchy-Dirichlet problem \eqref{pdd}, we assume that the nonlinear tensor $a\colon \Omega_{T}\times \mathbb{R}^{n}\to \mathbb{R}^{n}$ satisfies:
\begin{flalign}\label{regg}
\begin{cases}
\ t\mapsto a(x,t,z)\quad &\mbox{measurable for all} \ \ x\in \Omega, z\in \mathbb{R}^{n}\\
\ x\mapsto a(x,t,z)\quad &\mbox{differentiable for all} \ \ t\in (0,T),z\in \mathbb{R}^{n}\\
\ z\mapsto a(x,t,z)\in C(\mathbb{R}^{n},\mathbb{R}^{n})\cap C^{1}(\mathbb{R}^{n}\setminus \{0\},\mathbb{R}^{n})\quad &\mbox{for all} \ \ (x,t)\in \Omega_{T}
\end{cases}
\end{flalign}
and
\begin{flalign}\label{ref}
\begin{cases}
\ \snr{a(x,t,z)}+(\mu^{2}+\snr{z}^{2})^{\frac{1}{2}}\snr{\partial_{z}a(x,t,z)}\le L\left[(\mu^{2}+\snr{z}^{2})^{\frac{p-1}{2}}+(\mu^{2}+\snr{z}^{2})^{\frac{q-1}{2}}\right]\\
\ \left[\partial_{z}a(x,t,z)\xi\cdot \xi\right]\ge \nu(\mu^{2}+\snr{z}^{2})^{\frac{p-2}{2}}\snr{\xi}^{2}\\
\ \snr{\partial_{x}a(x,t,z)}\le \gamma(x,t) \left[(\mu^{2}+\snr{z}^{2})^{\frac{p-1}{2}}+(\mu^{2}+\snr{z}^{2})^{\frac{q-1}{2}}\right], 
\end{cases}
\end{flalign}
which holds for all $(x,t)\in \Omega_{T}$ and $z,\xi\in \mathbb{R}^{n}$. In \eqref{ref}, $\mu\in [0,1]$ is any number, exponents $(p,q)$ are so that
\begin{flalign}\label{pq}
q<p+2\left(\frac{1}{n+2}-\frac{p}{2d}\right)\qquad\mbox{with}\qquad p> \frac{2nd}{(n+2)(d-2)}
\end{flalign}
and 
\begin{flalign}\label{gamma}
\gamma\in L^{d}(\Omega_{T}) \ \ \mbox{for some} \ \ d>\max\left\{\frac{p}{2},1\right\}(n+2).
\end{flalign}
Finally, the function $f\colon \mathbb{R}^{n}\times \mathbb{R}\to \mathbb{R}$ satisfies
\begin{flalign}\label{ggg}
f\in C_{\loc}(\mathbb{R};L^{2}_{\loc}(\mathbb{R}^{n}))\cap L^{r}_{\loc}(\mathbb{R};W^{1,r}_{\loc}(\mathbb{R}^{n})),\quad \partial_{t}f\in L^{p'}_{\loc}(\mathbb{R};W^{-1,p}_{\loc}(\mathbb{R}^{n})),
\end{flalign}
where $r:=p'(q-1)$. In this setting, we define a weak solution to \eqref{pdd} as follows.
\begin{definition}\label{d.1}
A function $u\in f+L^{p}(0,T;W^{1,p}_{0}(\Omega))$ is a weak solution of problem \eqref{pdd} if and only if the identity
\begin{flalign}\label{wfp}
\int_{\Omega_{T}}\left[u\partial_{t}\varphi-a(x,t,Du)\cdot D\varphi\right] \ \dy=0
\end{flalign}
holds true for all $\varphi\in C^{\infty}_{c}(\Omega_{T})$ and, in addition, $u(\cdot,0)=f(\cdot,0)$ in the $L^{2}$-sense, i.e.:
\begin{flalign}\label{bf}
\lim_{\delta\to 0}\frac{1}{\delta}\int_{0}^{\delta}\int_{\Omega}\snr{u(x,s)-f(x,0)}^{2} \ \dx\ds=0.
\end{flalign}
\end{definition}
\begin{remark}
\emph{Let us compare the bound in \eqref{pq} with the one in force in the elliptic setting, i.e.:
\begin{flalign}\label{epq}
q<p+p\left(\frac{1}{n}-\frac{1}{d}\right),
\end{flalign}
see \cite{dm1,dm,elemarmas}. The restriction imposed in \eqref{pq} looks the right one: in fact, due to the different scaling in time, in \eqref{epq} $n$ must be replaced by $n+2$. Moreover, the usual parabolic deficit coming from the growth of the diffusive part of the equation affects also $d$:}
\begin{flalign*}
q<p+p\left(\frac{1}{n+2}-\left(d\cdot \frac{2}{p}\right)^{-1}\right)\cdot \frac{2}{p}.
\end{flalign*}
\emph{If we let $d\to \infty$ in \eqref{pq} and reverse the transformation prescribed by the caloric deficit phenomenon, we obtain}
\begin{flalign*}
q<p+\frac{p}{n},
\end{flalign*}
\emph{which is the same appearing in \cite{sharp} when the space-depending coefficient is Lipschitz-continuous.}
\end{remark}

\subsection{Auxiliary results}
In this section we collect some well-known facts that will have an important role throughout the paper.
\subsubsection*{On Sobolev functions} Let $w\in L^{1}(\Omega_{T},\mathbb{R}^{k})$, $k\ge 1$ be any function. If $h \in \mathbb{R}^{n}$ is a vector, we denote by $\tau_{h}\colon L^{1}(\Omega_{T},\mathbb{R}^{k})\to L^{1}(\Omega_{\snr{h}}\times (0,T),\mathbb{R}^{k})$ the standard finite difference operator in space, pointwise defined as
\begin{flalign*}
\tau_{h}w(x):=w(x+h,t)-w(x,t) \quad \mbox{for a.e.} \ (x,t)\in \Omega_{\snr{h}}\times(0,T),
\end{flalign*}
where $\Omega_{|h|}:=\{x \in \Omega \, : \, 
\dist(x, \partial \Omega) > |h|\}$ and by $\Di\colon L^1(\Omega_{T},\mathbb{\er}^{k}) \to L^{1}(\Omega_{|h|}\times (0,T),\mathbb{R}^{k})$ the spacial difference quotient operator, i.e.:
\begin{flalign*}
\Di w(x,t):=\frac{w(x+h,t)-w(x,t)}{\snr{h}}=\snr{h}^{-1}(\tau_{h}w(x,t)).
\end{flalign*}
Moreover, if $\tilde{h}\in \mathbb{R}$ is a number so that $\snr{h}<T$, we also recall the definition of finite difference operator in time $\tilde{\tau}_{\tilde{h}}\colon L^{1}(\Omega_{T})\to L^{1}(\Omega\times (\snr{\tilde{h}},T-\snr{\tilde{h}}))$:
\begin{flalign*}
\tilde{\tau}_{\tilde{h}}w(x,t):=w(x,t+h)-w(x,t).
\end{flalign*}
An important property of translation operators is their continuity in Lebesgue spaces.
\begin{lemma}\label{transc}
Let $\varphi\in C^{\infty}_{c}(\Omega)$ be any map, $h\in \mathbb{R}^{n}$ so that $\snr{h}\in \left(0,\frac{\dist(\supp(\varphi),\partial \Omega)}{4}\right)$ and $w\in L^{s}_{\loc}(\Omega_{T},\mathbb{R}^{k})$ with $s\in [1,\infty)$ and $k\in \N$. Then
\begin{flalign*}
\nr{(w(\ \cdot\ +h,t)-w(\cdot,t))\varphi}_{L^{s}(\Omega)}\to_{\snr{h}\to 0}0.
\end{flalign*}
\end{lemma}
It is also useful to recall a basic property of difference quotient.
\begin{lemma}\label{diffquo}
Let $w\in L^{1}_{\loc}(\Omega_{T})$ be any function. There holds that
\begin{itemize}
    \item if $w\in L^{s}_{\loc}(0,T;W^{1,s}_{\loc}(\Omega,\mathbb{R}^{k}))$, $s\in [1,\infty)$ and $\tilde{\Omega}\Subset$ is any open set, then
    \begin{flalign*}
    \nr{\Di w(\cdot,t)-Dw(\cdot,t)}_{L^{s}(\tilde{\Omega})}\to_{\snr{h}\to 0};
    \end{flalign*}
    \item if in addition $s>1$ and $\tilde{\Omega}\Subset \Omega$ is any open set so that 
    $$\sup_{\snr{h}>0}\int_{0}^{T}\int_{\tilde{\Omega}}\snr{\Di w(x,t)}^{s} \ \dx\dtt<\infty,$$
    then $Dw\in L^{s}(\ti{\Omega}\times (0,T))$ and $\nr{\Di w(\cdot,t)-Dw(\cdot,t)}_{L^{s}(\tilde{\Omega})}\to_{\snr{h}\to 0}0$.
\end{itemize}
\end{lemma}
When dealing with parabolic PDE, solutions in general posses a modest degree of regularity in the time-variable, and, in particular, time derivatives exist only in the distributional sense. For this reason, we recall the definition and main properties of Steklov averages, see e.g. \cite[Chapter 1]{d}.
\begin{definition}
Let $w\in L^{1}(\Omega_{T},\mathbb{R}^{k})$, $k\in \N$, be any function. For $\delta\in (0,T)$, the Steklov averages of $w$ are defined as
\begin{flalign*}
w_{\delta}:=\begin{cases}
\ \frac{1}{\delta}\int_{t}^{t+\delta}w(x,s) \ \ds\quad &t\in (0,T-\delta]\\
\ 0\quad &t>T-\delta
\end{cases}\quad \mbox{and}\quad w_{\bar{\delta}}:=\begin{cases}
\ \frac{1}{\delta}\int_{t-\delta}^{t}w(x,s) \ \ds\quad &t\in (\delta,T]\\
\ 0\quad &t<\delta.
\end{cases}
\end{flalign*}
\end{definition}
\begin{lemma}
If $w\in L^{s}_{\loc}(\Omega_{T})$, then $w_{\delta}\to_{\delta\to 0}w$ in $L^{s}_{\loc}(\Omega_{T-\varepsilon})$ for all $\varepsilon\in (0,T)$. If $w\in C(0,T; L^{s}(\Omega))$, then as $\delta\to 0$, $w_{\delta}(\cdot,t)$ converges to $w(\cdot,t)$ for all $t\in (0,T-\varepsilon)$ and all $\varepsilon\in (0,T)$. A similar statement holds for $w_{\bar{\delta}}$ as well.
\end{lemma}
We also record the definition of fractional Sobolev spaces.
\begin{definition}
A function $w\in L^{s}(\Omega_{T},\mathbb{R}^{k})$ belongs to the fractional Sobolev space $W^{\alpha,\theta;s}(\Omega_{T},\mathbb{R}^{k})$, $\alpha,\theta\in (0,1)$, $k\in \N$ provided that
\begin{flalign*}
\int_{0}^{T}\int_{\Omega}\int_{\Omega}\frac{\snr{w(x_{1},t)-w(x_{2},t)}^{s}}{\snr{x_{1}-x_{2}}^{n+s\alpha}} \ \dx_{1}\dx_{2}\dtt+\int_{0}^{T}\int_{0}^{T}\int_{\Omega}\frac{\snr{w(x,t_{1})-w(x,t_{2})}^{s}}{\snr{t_{1}-t_{2}}^{1+s\theta}} \ \dx\dtt_{1}\dtt_{2}<\infty.
\end{flalign*}
The local variant of $W^{\alpha,\theta;s}(\Omega_{T},\mathbb{R}^{k})$ can be defined in the usual way.
\end{definition}
The usual relation between Nikolski spaces and fractional Sobolev spaces holds in the parabolic setting as well.
\begin{proposition}\label{fracsob}
Let $w\in L^{s}(\Omega_{T},\mathbb{R}^{k})$, $(t_{1},t_{2})\Subset (0,T)$, $\tilde{\Omega}\Subset \Omega$ be an open set, $h\in \mathbb{R}^{n}$ be any vector with $\snr{h}<\frac{\dist(\tilde{\Omega},\partial \Omega)}{4}$ and $\tilde{h}\in \mathbb{R}$ be a number so that $\snr{\ti{h}}<\frac{\min\{t_{1},T-t_{2}\}}{4}$. Assume that
\begin{flalign*}
\int_{t_{1}}^{t_{2}}\int_{\tilde{\Omega}}\snr{w(x,t+\tilde{h})-w(x,t)} \ \dx\dtt\le c'\snr{\tilde{h}}^{s\theta}\quad \mbox{for some}  \ \ \theta \in (0,1),
\end{flalign*}
where $c'$ is a positive, absolute constant. Then there exists a constant $\tilde{c}\equiv \tilde{c}(n,s,c',\iota,t_{1},T-t_{2})>0$ such that
\begin{flalign*}
\int_{t_{1}}^{t_{2}}\int_{t_{1}}^{t_{2}}\int_{\tilde{\Omega}}\frac{\snr{w(x,l_{1})-w(x,l_{2})}^{s}}{\snr{l_{1}-l_{2}}^{1+s\iota}} \ \dx\dll_{1}\dll_{2}\le \tilde{c}<\infty\quad \mbox{for all} \ \ \iota \in (0,\theta).
\end{flalign*}
Suppose that
\begin{flalign*}
\int_{t_{1}}^{t_{2}}\int_{\tilde{\Omega}}\snr{w(x+h,t)-w(x,t)}^{s} \ \dx\dtt \le c'\snr{h}^{s\alpha}\quad \mbox{for some} \ \ \alpha\in (0,1),
\end{flalign*}
with $c'$ positive, absolute constant. Then,
\begin{flalign*}
\int_{t_{1}}^{t_{2}}\int_{\tilde{\Omega}}\int_{\tilde{\Omega}}\frac{\snr{w(x_{1},t)-w(x_{2},t)}^{s}}{\snr{x_{1}-x_{2}}^{n+s\gamma}} \ \dx_{1}\dx_{2}\dtt\le \tilde{c}<\infty\quad \mbox{for all} \ \ \gamma\in (0,\alpha),
\end{flalign*}
with $\tilde{c}\equiv \tilde{c}(n,s,c',\gamma,\dist(\ti{\Omega},\partial \Omega))$.
\end{proposition}
We refer to \cite{akm,pala,dumi1,lsu} for more details on this matter. We close this part with a fundamental compactness criterion in parabolic Sobolev spaces, whose proof can be found in \cite{sim}.
\begin{lemma}\label{al}
Let $X\subset B\subset Y$ be three Banach spaces such that the immersion $X\hookrightarrow B$ is compact and $1\le a_{1}\le a_{2}\le \infty$ be numbers satisfying the balance condition $a_{1}>a_{2}/(1+\sigma a_{2})$ for some $\sigma\in (0,1)$. If the set $\mathcal{J}$ is bounded in $L^{a_{2}}(0,T;X)\cap W^{\sigma,a_{1}}(0,T;Y)$, then $\mathcal{J}$ is compact in $L^{a_{2}}(0,T;B)$ and eventually in $C(0,T;B)$ when $a_{2}=\infty$.
\end{lemma}
\subsubsection*{Tools for $p$-laplacean type problems} For a constant $\tilde{c}\in [0,1]$ and $z\in \mathbb{R}^{n}$ we introduce the auxiliary vector field
\begin{flalign*}
V_{\tilde{c},s}(z):=(\tilde{c}^{2}+\snr{z}^{2})^{\frac{s-2}{4}}z\qquad s\in \{p,q\},
\end{flalign*}
which turns out to be very convenient in handling the monotonicity properties of certain operators.
\begin{lemma}\label{l1}
For any given $z_{1},z_{2}\in \mathbb{R}^{n}$, $z_{1}\not=z_{2}$ there holds that
\begin{flalign*}
\snr{V_{\tilde{c},s}(z_{1})-V_{\tilde{c},s}(z_{2})}^{2}\sim (\tilde{c}^{2}+\snr{z_{1}}^{2}+\snr{z_{2}}^{2})^{\frac{s-2}{2}}\snr{z_{1}-z_{2}}^{2},
\end{flalign*}
where the constants implicit in "$\sim$" depend only from $(n,s)$.
\end{lemma}
Another useful result is the following
\begin{lemma}\label{l6}
Let $s>-1$, $\tilde{c}\in [0,1]$ and $z_{1},z_{2}\in \mathbb{R}^{n}$ be so that $\tilde{c}+\snr{z_{1}}+\snr{z_{2}}>0$. Then
\begin{flalign*}
\int_{0}^{1}\left[\tilde{c}^{2}+\snr{z_{1}+\lambda(z_{2}-z_{1})}^{2}\right]^{\frac{s}{2}} \ \d\lambda\sim (\tilde{c}^{2}+\snr{z_{1}}^{2}+\snr{z_{2}}^{2})^{\frac{s}{2}},
\end{flalign*}
with constants implicit in "$\sim$" depending only from $s$. 
\end{lemma}
Finally, the iteration lemma.
\begin{lemma}\label{l0}
Let $\mathcal{Z}\colon [\rr,R)\to [0,\infty)$ be a function which is bounded on every interval $[\varrho, R_*]$ with $R_*<R$. Let $\varepsilon\in (0,1)$, $a_1,a_2,\gamma_{1},\gamma_{2}\ge 0$ be numbers. If
\begin{flalign*}
\mathcal{Z}(\tau_1)\le \varepsilon \mathcal{Z}(\tau_2)+ \frac{a_1}{(\tau_2-\tau_1)^{\gamma_{1}}}+\frac{a_2}{(\tau_2-\tau_1)^{\gamma_{2}}}\ \ \mbox{for all} \ \rr\le \tau_1<\tau_2< R\;,
\end{flalign*}
then
\begin{flalign*}
\mathcal{Z}(\rr)\le c\left[\frac{a_1}{(R-\rr)^{\gamma_{1}}}+\frac{a_2}{(R-\rr)^{\gamma_{2}}}\right]\;,
\end{flalign*}
holds with $c\equiv c(\varepsilon,\gamma_{1},\gamma_{2})$.
\end{lemma}

\section{Higher Sobolev regularity for non-degenerate systems}\label{high}
In this section we prove the existence of a suitably regular weak solution to Cauchy-Dirichlet problem 
\begin{flalign}\label{pdd+}
\begin{cases}
\ \partial_{t}v-\diver\ \tilde{a}(x,t,Dv)=0\quad &\mbox{in} \ \ \Omega_{T}\\
\ v=f\quad &\mbox{on} \ \ \partial_{par}\Omega_{T},
\end{cases}
\end{flalign}
where $f$ is as in \eqref{ggg} and the diffusive tensor $\tilde{a}\colon \Omega_{T}\times \mathbb{R}^{n}\to \mathbb{R}$ satisfies
\begin{flalign}\label{regg+}
\begin{cases}
\ t\mapsto \tilde{a}(x,t,z)\quad &\mbox{measurable for all} \ \ x\in \Omega, z\in \mathbb{R}^{n}\\
\ x\mapsto \tilde{a}(x,t,z)\quad &\mbox{differentiable for all} \ \ t\in (0,T),z\in \mathbb{R}^{n}\\
\ z\mapsto \tilde{a}(x,t,z)\in C^{1}(\mathbb{R}^{n},\mathbb{R}^{n})\quad &\mbox{for all} \ \ (x,t)\in \Omega_{T}
\end{cases}
\end{flalign}
and
\begin{flalign}\label{ref+}
\begin{cases}
\ \snr{\tilde{a}(x,t,z)}+(\ti{\mu}^{2}+\snr{z}^{2})^{\frac{1}{2}}\snr{\partial_{z}\ti{a}(x,t,z)}\le L\left[(\ti{\mu}^{2}+\snr{z}^{2})^{\frac{p-1}{2}}+(\ti{\mu}^{2}+\snr{z}^{2})^{\frac{q-1}{2}}\right]\\
\ \left[\partial_{z}\ti{a}(x,t,z)\xi\cdot \xi\right]\ge \nu(\ti{\mu}^{2}+\snr{z}^{2})^{\frac{p-2}{2}}\snr{\xi}^{2}\\
\ \snr{\partial_{x}\ti{a}(x,t,z)}\le \gamma(x,t) \left[(\ti{\mu}^{2}+\snr{z}^{2})^{\frac{p-1}{2}}+(\ti{\mu}^{2}+\snr{z}^{2})^{\frac{q-1}{2}}\right], 
\end{cases}
\end{flalign}
for all $(x,t)\in \Omega_{T}$ and $z,\xi\in \mathbb{R}^{n}$. In \eqref{ref+}, $(p,q)$ are linked by the relation in \eqref{pq}, $\gamma$ is as in \eqref{gamma} and 
\begin{flalign}\label{extra}
\ti{\mu}>0.
\end{flalign}
Our first result is the following
\begin{proposition}\label{p1}
Let $f\colon \mathbb{R}^{n}\times \mathbb{R}\to \mathbb{R}$ be as in \eqref{ggg} and $\ti{a}\colon \Omega_{T}\times \mathbb{R}^{n}\to \mathbb{R}^{n}$ be a Carath\'eodory vector field satisfying \eqref{regg+}, \eqref{ref+}, \eqref{pq}, \eqref{gamma} and \eqref{extra}. Then there exists a weak solution $v\in L^{p}(0,T;W^{1,p}(\Omega))$ of Cauchy-Dirichlet problem \eqref{pdd+} such that
\begin{flalign}\label{sss}
v\in L^{s}_{\loc}(0,T;W^{1,s}_{\loc}(\Omega))\quad \mbox{for all} \ \ s\in \left[1,p+\frac{4}{\tilde{n}}\right]
\end{flalign}
satisfying
\begin{flalign}\label{cv7}
\partial_{t}v\in L^{l}_{\loc}(\Omega_{T})\quad \mbox{for some} \ \ l\equiv l(n,p,q,d)\in \left(1,\min\{2,p\}\right)
\end{flalign}
and
\begin{flalign}\label{difff}
Dv\in L^{\infty}_{\loc}(0,T,L^{2}_{\loc}(\Omega,\mathbb{R}^{n}))\quad \mbox{with}\quad V_{p}(Dv)\in L^{2}_{\loc}(0,T;W^{1,2}_{\loc}(\Omega,\mathbb{R}^{n})).
\end{flalign}
\end{proposition}
For the sake of simplicity, we shall split the proof of Proposition \ref{p1} into eight steps.
\subsubsection*{Step 1: Approximating Cauchy-Dirichlet problems}
For the ease of notation, we define numbers:
\begin{flalign}\label{r}
m:=\frac{d}{d-2}>1,\qquad \tilde{q}:=\max\left\{q-\frac{p}{2},1\right\},
\end{flalign}
and, for $j\in \N$, consider a usual family of non-negative mollifiers $\{\psi_{j}\}$ of $\mathbb{R}^{n+1}$. We then regularize $f$ via convolution against $\{\psi_{j}\}$, thus obtaining the sequence $\{f_{j}\}:=\{f*\psi_{j}\}$, set 
\begin{flalign}\label{ej}
\varepsilon_{j}:=\left(1+j+\nr{f_{j}}_{L^{2m\tilde{q}}(\Omega_{T})}^{2m\tilde{q}}\right)^{-1},\qquad \ti{H}(z):=(\ti{\mu}^{2}+\snr{z}^{2}),
\end{flalign}
correct the nonstandard growth of the diffusive tensor $\ti{a}(\cdot)$ as follows:
\begin{flalign}\label{defa}
\ti{a}_{j}(x,t,z):=\ti{a}(x,t,z)+\varepsilon_{j}\ti{H}(z)^{\frac{2m\tilde{q}-2}{2}}z
\end{flalign}
and consider solutions $v_{j}\in L^{2m\ti{q}}(0,T;W^{1,2m\ti{q}}(\Omega))$ of the following Cauchy-Dirichlet problem
\begin{flalign}\label{pd}
\begin{cases}
\ \partial v_{j}-\diver \ \taa_{j}(x,t,Dv_{j})=0\quad &\mbox{in} \ \ \Omega_{T}\\
\ v_{j}=f_{j}\quad &\mbox{on} \ \ \partial_{par}\Omega_{T}.
\end{cases}
\end{flalign}
By \eqref{extra}, \eqref{ref} and the definition in \eqref{defa}, it can be easily seen that \eqref{regg+} holds and
and
\begin{flalign}\label{refreg}
\begin{cases}
\ \snr{\taa_{j}(x,t,z)}+\ti{H}(z)^{\frac{1}{2}}\snr{\partial_{z}\taa_{j}(x,t,z)}\le c\left[\ti{H}(z)^{\frac{p-1}{2}}+\ti{H}(z)^{\frac{q-1}{2}}\right]+c\varepsilon_{j}\ti{H}(z)^{\frac{2m\tilde{q}-1}{2}}\\
\ \partial_{z}\taa_{j}(x,t,z)\ge c\left[\ti{H}(z)^{\frac{p-2}{2}}\varepsilon_{j}\ti{H}(z)^{\frac{2m\tilde{q}-2}{2}}\right]\snr{\xi}^{2}\\
\ \snr{\partial_{x}\taa_{j}(x,t,z)}\le c\gamma(x,t)\left[\ti{H}(z)^{\frac{p-1}{2}}+\ti{H}(z)^{\frac{q-1}{2}}\right],
\end{cases}
\end{flalign}
for all $(x,t)\in \Omega_{T}$, $z,\xi\in \mathbb{R}^{n}$, with $\gamma$ as in \eqref{gamma} and $c\equiv c(n,\nu,L,p,q,d)$. We recall that the weak formulation associated to problem \eqref{pdd+} reads as
\begin{flalign}\label{wfj}
\int_{\Omega_{T}}\left[v_{j}\partial_{t}\varphi -\taa(x,t,Dv_{j})\cdot D\varphi\right] \ \dy=0\quad \mbox{for all} \ \ \varphi\in C^{\infty}_{c}(\Omega_{T})
\end{flalign}
and the attainment of the boundary datum $f_{j}$ must be considered in the $L^{2}$-sense as in Definition \ref{d.1}.
\subsubsection*{Step 2: Uniform energy bounds} Our main goal it to prove that the sequence $\{v_{j}\}$ is bounded, uniformly with respect to $j\in \N$ in the space-time $L^{p}$-norm. Since this is quite a routine procedure, we will just sketch it and refer the reader to \cite{bdm,si}, for more details. Modulo using Steklov averages, we can test \eqref{wfj} against the difference $v_{j}-f_{j}$ to get
\begin{flalign}\label{24}
\int_{\Omega}&\snr{v_{j}(x,t)-f_{j}(x,t)}^{2} \ \dx\nonumber \\
&+\int_{0}^{t}\int_{\Omega}\taa_{j}(x,s,Dv_{j})\cdot(Dv_{j}-Df_{j}) \ \dx \ds \nonumber \\
=&-\int_{0}^{t}\langle v_{j}-f_{j},\partial_{t}f_{j}\rangle_{W^{1,p}_{0}(\Omega)} \ \ds\quad \mbox{for a.e.} \ \ t\in (0,T).
\end{flalign}
By $\eqref{refreg}_{2}$, H\"older and Young inequalities, if $p\ge 2$ a straightforward computation renders that 
\begin{flalign*}
\int_{0}^{t}\int_{\Omega}&\snr{Dv_{j}}^{p} \ \dx\ds +\varepsilon_{j}\int_{0}^{t}\int_{\Omega}\snr{Dv_{j}}^{2m\tilde{q}} \ \dx\ds\nonumber \\
\lesssim &\int_{0}^{t}\int_{\Omega}\left[\taa_{j}(x,s,Dv_{j})-\taa_{j}(x,s,Df_{j})\right](Dv_{j}-Df_{j}) \ \dx\ds\nonumber \\
&+\int_{0}^{t}\int_{\Omega}\left[\snr{Df_{j}}^{p}+\varepsilon_{j}\snr{Df_{j}}^{2m\ti{q}}\right] \ \dx\ds,
\end{flalign*}
while if $1<p<2$ there holds that
\begin{flalign*}
\int_{0}^{t}\int_{\Omega}&\snr{Dv_{j}}^{p} \ \dx\ds+\varepsilon_{j}\int_{0}^{t}\int_{\Omega}\snr{Dv_{j}}^{2m\ti{q}} \ \dx\ds\nonumber \\
\lesssim &\frac{1}{\sigma}\int_{0}^{t}\int_{\Omega}\left[\taa_{j}(x,s,Dv_{j})-a_{j}(x,s,Df_{j})\right]\cdot(Dv_{j}-Df_{j}) \ \dx\ds\nonumber \\
&+\sigma \int_{0}^{t}\int_{\Omega}\snr{Dv_{j}}^{p} \ \dy+\int_{\Omega_{t_{0}}}\left[\snr{Df_{j}}^{p}+\varepsilon_{j}\snr{Df_{j}}^{2m\ti{q}}\right] \ \dx\ds.
\end{flalign*}
Moreover, using $\eqref{refreg}_{1}$, H\"older and Young inequalities we have
\begin{flalign*}
\int_{0}^{t}\int_{\Omega}&\taa_{j}(x,s,Df_{j})\cdot(Dv_{j}-Df_{j}) \ \dx\ds\lesssim \sigma \int_{0}^{t}\int_{\Omega}\snr{Dv_{j}}^{p} \ \dx\ds+\sigma\varepsilon_{j}\int_{0}^{t}\int_{\Omega}\snr{Dv_{j}}^{2m\ti{q}} \ \dx\ds \nonumber\\
&+\frac{1}{\sigma}\int_{0}^{t}\int_{\Omega}\left[1+\snr{Df_{j}}^{r}\right] \ \dx\ds+\frac{\varepsilon_{j}}{\sigma}\int_{0}^{t}\int_{\Omega}\ti{H}(Df_{j})^{m\tilde{q}} \ \dx\ds.
\end{flalign*}
Here, we also used that $q\ge p\Rightarrow r\ge q$ and, of course, that $2m\ti{q}> 2$. Finally, by H\"older, Sobolev-Poincar\'e and Young inequalities
\begin{flalign*}
\left| \ \int_{0}^{t}\langle v_{j}-f_{j}\rangle_{W^{1,p}_{0}(\Omega)} \ \ds \ \right|\lesssim \sigma\int_{0}^{t}\int_{\Omega}\snr{Dv_{j}}^{p} \ \dx\ds+\frac{1}{\sigma}\nr{\partial_{t}f_{j}}^{p'}_{L^{p'}(0,t_{1};W^{-1,p'}(\Omega))}.
\end{flalign*}
Inserting the content of all the previous displays in \eqref{24}, recalling \eqref{ej}, \eqref{ggg} and well-known convolution properties, choosing $\sigma>0$ small enough, we obtain
\begin{flalign}\label{unibd1}
\int_{0}^{t}\int_{\Omega}&\snr{Dv_{j}}^{p} \ \dx\ds+\varepsilon_{j}\int_{0}^{t}\int_{\Omega}\snr{Dv_{j}}^{2m\ti{q}} \ \dx\ds+\int_{\Omega}\snr{v_{j}(x,t)-f_{j}(x,t)}^{2} \ \dx\nonumber \\
\lesssim& \left[\int_{\int_{0}^{t}\Omega}\left[1+\snr{Df_{j}}^{r}\right] \ \dx\ds+\varepsilon_{j}\int_{0}^{t}\int_{\Omega}\ti{H}(Df_{j})^{m\ti{q}} \ \dx\ds+\nr{\partial_{t}f_{j}}^{p'}_{L^{p'}(0,t;W^{-1,p'}(\Omega))}\right]\nonumber \\
\lesssim &\left[\int_{0}^{t}\int_{\Omega}\left[1+\snr{Df}^{r}\right] \ \dx\ds+\nr{\partial_{t}f}^{p'}_{L^{p'}(0,t;W^{-1,p'}(\Omega))}+1\right]\nonumber\\
\lesssim &\left[\nr{Df}_{L^{r}(\Omega_{T})}+\nr{\partial_{t}f}^{p'}_{L^{p'}(0,T;W^{-1,p'}(\Omega))}+1\right].
\end{flalign}
As stated at the end of Section \ref{ma}, none of the constants implicit in "$\lesssim$" depends on $t\in (0,T)$, therefore we can send $t\to T$ on the right-hand side of \eqref{unibd1} to get
\begin{flalign}\label{unibd}
&\int_{0}^{t}\int_{\Omega}\snr{Dv_{j}}^{p} \ \dx\ds+\varepsilon_{j}\int_{0}^{t}\int_{\Omega}\snr{Dv_{j}}^{2m\ti{q}} \ \dx\ds+\int_{\Omega}\snr{v_{j}(x,t)-f_{j}(x,t)}^{2} \ \dx\nonumber \\
&\qquad \lesssim \left[\nr{Df}_{L^{r}(\Omega_{T})}^{r}+\nr{\partial_{t}f}^{p'}_{L^{p'}(0,T;W^{-1,p'}(\Omega))}+1\right]=:\mathcal{C}_{f}.
\end{flalign}
\subsubsection*{Step 3: Caccioppoli inequality}
Let $h\in \mathbb{R}^{n}\setminus \{0\}$ any vector satisfying $\snr{h}\in (0,1)$, $B_{\rr}\subset \Omega$ a ball with radius $0<\rr\le 1$ and so that $B_{2\rr}\Subset \Omega$, $g\in W^{1,\infty}(\mathbb{R})$ a non-negative function with bounded, piecewise continuous, non-negative first derivative and $\chi\in W^{1,\infty}([0,T])$ with $\chi(0)=0$, $\varphi\in C^{\infty}(B_{\rr},[0,1])$ two cut-off functions. By the approximation procedure developed  e.g. in \cite[Section 3]{bdm} or \cite[Section 3.1]{s}, we can test \eqref{wfj} against a suitably regularized version of the comparison map $\varphi^{2}\chi\Di u_{j}g(\snr{\Di u_{j}}^{2})$ and manipulate it to obtain, for a.e. $\tau\in (0,\min\{T,1\})$,
\begin{flalign}\label{0}
\frac{1}{2}\int_{B_{\rr}}&\varphi^{2}\chi\left(\int_{0}^{\snr{\Di v_{j}}^{2}}g(s) \ \ds\right) \ \dx+\sum_{k=1}^{n}\int_{Q_{\tau}}\varphi^{2}\chi\Di\taa_{j}^{k}(x,t,Dv_{j})D_{k}\left[\Di v_{j}g(\snr{\Di v_{j}}^{2})\right] \ \dy\nonumber \\
=&-2\sum_{k=1}^{n}\int_{Q_{\tau}}\varphi\chi\left(\Di v_{j}g(\snr{\Di v_{j}}^{2})\right)\Di\taa_{j}^{k}(x,t,Dv_{j}) D_{k}\varphi \ \dy\nonumber \\
&-\frac{1}{2}\int_{Q_{\tau}}\left(\int_{0}^{\snr{\Di v_{j}}^{2}}g(s)\ \ds\right)\varphi^{2}\partial_{t}\chi \ \dy,
\end{flalign}
where we abbreviated $Q_{\tau}:=B_{\rr}\times (0,\tau)$. We also reduce further the size of $\snr{h}$: we ask that
\begin{flalign}\label{hhh}
\snr{h}\in \left(0,\frac{\dist(\supp(\varphi),\partial B_{\rr})}{10000}\right).
\end{flalign}
Using the mean value theorem, we rearrange $\Di\ta(x,t,Du_{j})$ in a more convenient way:
\begin{flalign*}
\Di\taa_{j}^{k}(x,t,Dv_{j})=&\snr{h}^{-1}\left[\tilde{a}^{k}(x+h,t,Dv_{j}(x+h))-\tilde{a}^{k}(x,t,Dv_{j}(x+h))\right]\nonumber \\
&+\snr{h}^{-1}\varepsilon_{j}\left[\ti{H}(Dv_{j}(x+h))^{\frac{2m\tilde{q}-2}{2}}D_{k}v_{j}(x+h)-\ti{H}(Dv_{j}(x))^{\frac{2m\tilde{q}-2}{2}}D_{k}v_{j}(x)\right]\nonumber \\
&+\snr{h}^{-1}\left[\tilde{a}^{k}(x,t,Dv_{j}(x+h))-\tilde{a}^{k}(x,t,Dv_{j}(x))\right]\nonumber \\
=&\snr{h}^{-1}\sum_{l=1}^{n}\left[\int_{0}^{1}\partial_{x_{l}}\tilde{a}^{k}(x+\lambda h,t,Dv_{j}(x+h))h^{l} \ \d \lambda\right]\nonumber \\
&+\sum_{l=1}^{n}\left[\int_{0}^{1}\partial_{z_{l}}\taa_{j}^{k}(x,t,Dv_{j}(x)+\lambda \tau_{h}Dv_{j}(x))\d\lambda\right]\Di D_{l}v_{j}.
\end{flalign*}
Plugging this expansion in \eqref{0} we eventually get
\begin{flalign}\label{1}
\frac{1}{2}\int_{B_{\rr}}&\varphi^{2}\chi\left(\int_{0}^{\snr{\Di v_{j}}^{2}}g(s) \ \ds\right) \ \dx\nonumber \\
&+\snr{h}^{-1}\sum_{k,l=1}^{n}\int_{Q_{\tau}}\varphi^{2}\chi\left[\int_{0}^{1}\partial_{x_{l}}\ti{a}^{k}(x+\lambda h,t,Dv_{j}(x+he_{i}))h^{l} \ \d\lambda\right]D_{k}\left[\Di v_{j}g(\snr{\Di v_{j}}^{2})\right] \ \dy\nonumber \\
&+\sum_{k,l=1}^{n}\int_{Q_{\tau}}\varphi^{2}\chi\left[\int_{0}^{1}\partial_{z_{l}}\taa_{j}^{k}(x,t,Dv_{j}(x)+\lambda \tau_{h}Dv_{j}(x)) \ \d \lambda\right]\Di D_{l}v_{j}D_{k}\left[\Di v_{j}g(\snr{\Di v_{j}}^{2})\right] \ \dy\nonumber \\
=&-2\snr{h}^{-1}\sum_{k,l=1}^{n}\int_{Q_{\tau}}\varphi\chi\left(\Di v_{j}g(\snr{\Di v_{j}}^{2})\right)\left[\int_{0}^{1}\partial_{x_{l}}\ti{a}^{k}(x+\lambda h,t,Dv_{j}(x+he_{i}))h^{l} \ \d\lambda\right]D_{k}\varphi \ \dy\nonumber \\
&-2\sum_{k,l=1}^{n}\int_{Q_{\tau}}\varphi\chi\left(\Di v_{j}g(\snr{\Di v_{j}}^{2})\right)\left[\int_{0}^{1}\partial_{z_{l}}\taa_{j}^{k}(x,t,Dv_{j}(x)+\lambda \tau_{h}v_{j}) \ \d\lambda\right]\Di D_{l}v_{j}D_{k}\varphi \ \dy\nonumber \\
&+\frac{1}{2}\int_{Q_{\tau}}\left(\int_{0}^{\snr{\Di v_{j}}^{2}}g(s) \ \ds\right)\varphi^{2}\partial_{t}\chi \ \dy.
\end{flalign}
For reasons that will be clear in a few lines, we introduce the shorthands
\begin{flalign*}
\mathcal{D}(h):=\left(\ti{\mu}^{2}+\snr{Dv_{j}(x+h)}^{2}+\snr{Dv_{j}(x)}^{2}\right)\quad \mbox{and}\quad \mathcal{G}(h):=\left(g(\snr{\Di v_{j}}^{2})+\snr{\Di v_{j}}^{2}g'(\snr{\Di v_{j}}^{2})\right),
\end{flalign*}
and notice that, by \eqref{extra}, $\mathcal{D}(h)>\ti{\mu}^{2}>0$. Now we start estimating all the terms appearing in \eqref{1}. For the sake of clarity, we split
\begin{flalign*}
&\mbox{(I)}:=\snr{h}^{-1}\sum_{k=1}^{n}\int_{Q_{\tau}}\varphi^{2}\chi\left[\int_{0}^{1}\partial_{x_{l}}\ti{a}^{k}(x+\lambda h,t,Dv_{j}(x+he_{i}))h^{l} \ \d\lambda\right]D_{k}\left[\Di v_{j}g(\snr{\Di v_{j}}^{2})\right] \ \dy\nonumber \\
&\ \ =\snr{h}^{-1}\sum_{k=1}^{n}\int_{Q_{\tau}}\varphi^{2}\chi\left[\int_{0}^{1}\partial_{x_{l}}\ti{a}^{k}(x+\lambda h,t,Dv_{j}(x+he_{i}))h^{l} \ \d\lambda\right]\Di D_{k}v_{j}g(\snr{\Di v_{j}}^{2}) \ \dy\nonumber \\
& \ \ +2\snr{h}^{-1}\sum_{k=1}^{n}\int_{Q_{\tau}}\varphi^{2}\chi\left[\int_{0}^{1}\partial_{x_{l}}\ti{a}^{k}(x+\lambda h,t,Dv_{j}(x+he_{i}))h^{l} \ \d\lambda\right]\snr{\Di v_{j} }^{2}g'(\snr{\Di v_{j}}^{2})\Di D_{k}v_{j} \ \dy\nonumber \\
& \ \ =:\mbox{(I)}_{1}+\mbox{(I)}_{2}.
\end{flalign*}
With \eqref{refreg}$_{3}$, \eqref{gamma}, H\"older and Young inequalities we bound
\begin{flalign*}
\snr{\mbox{(I)}_{1}}&+\snr{\mbox{(I)}_{2}}\le c\int_{Q_{\tau}}\varphi^{2}\chi\left(\int_{0}^{1}\gamma(x+\lambda h,t) \ \d\lambda\right)\left[\mathcal{D}(h)^{\frac{p-1}{2}}+\mathcal{D}(h)^{\frac{q-1}{2}}\right]\mathcal{G}(h)\snr{\Di Dv_{j}} \ \dy\nonumber \\
\le&\sigma \int_{Q_{\tau}}\varphi^{2}\chi\mathcal{G}(h)\mathcal{D}(h)^{\frac{p-2}{2}}\snr{\Di Dv_{j}}^{2} \ \dy\nonumber \\
&+\frac{c}{\sigma}\int_{Q_{\tau}}\varphi^{2}\chi\left(\int_{0}^{1}\gamma(x+\lambda h,t) \ \d\lambda\right)^{2}\left[\mathcal{D}(h)^{\frac{p}{2}}+\mathcal{D}(h)^{\frac{2q-p}{2}}\right]\mathcal{G}(h) \ \dy\nonumber \\
\le&\sigma \int_{Q_{\tau}}\varphi^{2}\chi\mathcal{G}(h)\mathcal{D}(h)^{\frac{p-2}{2}}\snr{\Di Dv_{j}}^{2} \ \dy\nonumber \\
&+\frac{c}{\sigma}\int_{0}^{\tau}\nr{\gamma(\cdot,t)}^{2}_{L^{d}(B_{2\rr})}\left(\int_{B_{\rr}}\varphi^{2m}\chi^{m}\left[\mathcal{D}(h)^{\frac{pm}{2}}+\mathcal{D}(h)^{\frac{(2q-p)m}{2}}\right]\mathcal{G}(h)^{m} \ \dx \right)^{\frac{1}{m}} \ \dt\nonumber \\
\le &\sigma \int_{Q_{\tau}}\varphi^{2}\chi\mathcal{G}(h)\mathcal{D}(h)^{\frac{p-2}{2}}\snr{\Di Dv_{j}}^{2} \ \dy+\frac{c}{\sigma}\left(\int_{Q_{\tau}}\varphi^{2m}\chi^{m}\left[1+\mathcal{D}(h)^{m\left(q-\frac{p}{2}\right)}\right]\mathcal{G}(h)^{m} \ \dy\right)^{\frac{1}{m}},
\end{flalign*}
for $c\equiv c(\texttt{data})$. Moreover, by \eqref{refreg}$_{2}$, Lemmas \ref{l1} and \ref{l6}, we obtain
\begin{flalign*}
&\mbox{(II)}:=\sum_{k,l=1}^{n}\int_{Q_{\tau}}\varphi^{2}\chi\left[\int_{0}^{1}\partial_{z_{l}}\taa_{j}^{k}(x,t,Dv_{j}(x)+\lambda \tau_{h}Dv_{j}(x)) \ \d \lambda\right]\Di D_{l}v_{j}D_{k}\left[\Di v_{j}g(\snr{\Di v_{j}}^{2})\right] \ \dy\nonumber \\
&\ \ =\sum_{k,l=1}^{n}\int_{Q_{\tau}}\varphi^{2}\chi\left[\int_{0}^{1}\partial_{z_{l}}\taa_{j}^{k}(x,t,Dv_{j}(x)+\lambda \tau_{h}Dv_{j}(x)) \ \d \lambda\right]\Di D_{l}v_{j}\Di D_{k} v_{j} g(\snr{\Di v_{j}}^{2}) \ \dy\nonumber \\
&\ \ +2\sum_{k,l=1}^{n}\int_{Q_{\tau}}\varphi^{2}\chi\left[\int_{0}^{1}\partial_{z_{l}}\taa_{j}^{k}(x,t,Dv_{j}(x)+\lambda \tau_{h}Dv_{j}(x)) \ \d \lambda\right]\Di D_{l}v_{j}\snr{\Di v_{j}}^{2}g'(\snr{\Di v_{j}}^{2})\Di D_{k}v_{j} \ \dy\nonumber \\
& \ \ \ge c\snr{h}^{-2}\int_{Q_{\tau}}\varphi^{2}\chi\mathcal{D}(h)^{\frac{p-2}{2}}\snr{\tau_{h} Dv_{j}}^{2}\mathcal{G}(h) \ \dy+c\snr{h}^{-2}\varepsilon_{j} \int_{Q_{\tau}}\varphi^{2}\chi\mathcal{D}(h)^{\frac{2m\tilde{q}-2}{2}}\snr{\tau_{h} Dv_{j}}^{2}\mathcal{G}(h) \ \dy\nonumber \\
& \ \ \ge c\int_{Q_{\tau}}\varphi^{2}\chi \mathcal{G}(h)\snr{\Di V_{\ti{\mu},p}(Dv_{j})}^{2} \ \dy+c\varepsilon_{j} \int_{Q_{\tau}}\varphi^{2}\chi\mathcal{G}(h)\snr{\Di V_{\ti{\mu},2m\tilde{q}}(Dv_{j})}^{2} \ \dy,
\end{flalign*}
with $c\equiv c(n,\nu,p,q,d)$. With \eqref{refreg}$_{1,3}$, H\"{o}lder and Young inequalities we finally obtain
\begin{flalign*}
&\snr{\mbox{(III)}}+\snr{\mbox{(IV)}}:=\nonumber \\
&\quad -2\snr{h}^{-1}\sum_{k,l=1}^{n}\int_{Q_{\tau}}\varphi\chi\left(\Di v_{j}g(\snr{\Di v_{j}}^{2})\right)\left[\int_{0}^{1}\partial_{x_{l}}\taa^{k}(x+h\lambda,t,Dv_{j}(x+j))h^{l} \ \d\lambda\right]D_{k}\varphi \ \dy\nonumber \\
&\quad-2\sum_{k,l=1}^{n}\int_{Q_{\tau}}\varphi\chi\left(\Di v_{j}g(\snr{\Di v_{j}}^{2})\right)\left[\int_{0}^{1}\partial_{z_{l}}\taa_{j}^{k}(x,t,Dv_{j}(x)+\lambda \tau_{h}Dv_{j}) \ \d\lambda\right]\Di D_{l}v_{j}D_{k}\varphi \ \dy\nonumber \\
&\quad\le \sigma\int_{Q_{\tau}}\varphi^{2}\chi\mathcal{G}(h)\mathcal{D}(h)^{\frac{p-2}{2}}\snr{\Di Dv_{j}}^{2} \ \dy+\sigma \varepsilon_{j}\int_{Q_{\tau}}\varphi^{2}\chi\mathcal{G}(h)\mathcal{D}(h)^{\frac{2m\tilde{q}-2}{2}}\snr{\Di Dv_{j}}^{2}\ \dy\nonumber \\
&\quad+\frac{c\varepsilon_{j}}{\sigma}\int_{Q_{\tau}}\chi\snr{D\varphi}^{2}g(\snr{\Di v_{j}}^{2})\snr{\Di v_{j}}^{2}\mathcal{D}(h)^{\frac{2m\tilde{q}-2}{2}} \ \dy\nonumber\\
&\quad +\frac{c}{\sigma}\int_{Q_{\tau}}\chi\snr{D\varphi}^{2}\snr{\Di v_{j}}^{2}g(\snr{\Di v_{j}}^{2})\left[\mathcal{D}(h)^{\frac{p-2}{2}}+\mathcal{D}(h)^{q-\frac{p}{2}-1}\right] \ \dy\nonumber \\
&\quad+c\nr{\gamma}^{2}_{L^{d}(\Omega_{T})}\left(\int_{Q_{\tau}}\chi^{m}\varphi^{2m}g(\snr{\Di v_{j}}^{2})^{m}\left[\mathcal{D}(h)^{\frac{pm}{2}}+\mathcal{D}(h)^{m\left(q-\frac{p}{2}\right)}\right] \ \dy\right)^{\frac{1}{m}},
\end{flalign*}
where $c\equiv c(\texttt{data})$. Merging the content of all the previous displays and choosing $\sigma>0$ sufficiently small, we end up with
\begin{flalign}\label{2}
\frac{1}{2}\int_{B_{\rr}}&\varphi^{2}\chi\left(\int_{0}^{\snr{\Di v_{j}}^{2}}g(s) \ \ds\right) \ \dx\nonumber \\
&+\int_{Q_{\tau}}\varphi^{2}\chi\mathcal{G}(h)\snr{\Di V_{\ti{\mu},p}(Dv_{j})}^{2} \ \dy+\varepsilon_{j}\int_{Q_{\tau}}\varphi^{2}\chi\mathcal{G}(h)\snr{\Di V_{\ti{\mu},2m\tilde{q}}(Dv_{j})}^{2} \ \dy\nonumber \\
\le &c\left(\int_{Q_{\tau}}\chi^{m}\varphi^{2m}\mathcal{G}(h)^{m}\left[1+\mathcal{D}(h)^{m\left(q-\frac{p}{2}\right)}\right] \ \dy\right)^{\frac{1}{m}}\nonumber \\
&+c\int_{Q_{\tau}}\chi\snr{D\varphi}^{2}\snr{\Di v_{j}}^{2}g(\snr{\Di v_{j}}^{2})\left[\mathcal{D}(h)^{\frac{p-2}{2}}+\mathcal{D}(h)^{q-\frac{p}{2}-1}\right] \ \dy\nonumber \\
&+c\varepsilon_{j}\int_{Q_{\tau}}\chi\snr{\Di v_{j}}^{2}g(\snr{\Di v_{j}}^{2})\mathcal{D}(h)^{\frac{2m\tilde{q}-2}{2}}\snr{D\varphi}^{2} \ \dy\nonumber \\
&+c\int_{Q_{\tau}}\left(\int_{0}^{\snr{\Di v_{j}}^{2}}g(s) \ \ds\right)\varphi^{2}\partial_{t}\chi \ \dy=:\mathcal{I}(h),
\end{flalign}
with $c\equiv c(\texttt{data})$. In \eqref{2}, we also used that $m>1$ and that, being $p\le q$ we have that $\frac{p}{2}\le \frac{q}{2}\le q-\frac{p}{2}$. For $z\in \mathbb{R}^{n}$, set $\hat{\mathcal{G}}(z):=\left(g(\snr{z}^{2})+\snr{z}^{2}g'(\snr{z}^{2})\right)$. Now we recall \eqref{hhh} and that $g(\cdot)$ is bounded with bounded, piecewise continuous, non-negative first derivative. Keeping \eqref{extra} in mind, it is then easy to see that by Lemmas \ref{transc}-\ref{diffquo}, we can use Fatou Lemma on the left-hand side of \eqref{2} and a well-known variant of the dominated convergence theorem on the right-hand side of \eqref{2} to end up with
\begin{flalign}\label{5}
&\frac{1}{2}\int_{B_{\rr}}\varphi^{2}\chi\left(\int_{0}^{\snr{D v_{j}}^{2}}g(s) \ \ds\right) \ \dx\nonumber \\
&\quad+\int_{Q_{\tau}}\varphi^{2}\chi\hat{\mathcal{G}}(Dv_{j})\snr{D V_{\ti{\mu},p}(Dv_{j})}^{2} \ \dy+\varepsilon_{j}\int_{Q_{\tau}}\varphi^{2}\chi\hat{\mathcal{G}}(Dv_{j})\snr{D V_{\ti{\mu},2m\tilde{q}}(Dv_{j})}^{2} \ \dy\nonumber \\
&\quad\le c\left(\int_{Q_{\tau}}\chi^{m}\left(\snr{D\varphi}^{2m}+\varphi^{2m}\right)\hat{\mathcal{G}}(Dv_{j})^{m}\left[1+\ti{H}(Dv_{j})^{m\left(q-\frac{p}{2}\right)}\right] \ \dy\right)^{\frac{1}{m}}\nonumber \\
&\quad+c\varepsilon_{j}\int_{Q_{\tau}}\chi\snr{D\varphi}^{2} g(\snr{Dv_{j}}^{2})\ti{H}(Dv_{j})^{m\tilde{q}} \ \dy\nonumber \\
&\quad+c\int_{Q_{\tau}}\left(\int_{0}^{\snr{Dv_{j}}^{2}}g(s) \ \ds\right) \ \varphi^{2}\partial_{t}\chi \ \dy,
\end{flalign}
with $c\equiv c(\texttt{data})$. 
\subsubsection*{Step 4: Higher weak differentiability and interpolation} Our starting point is inequality \eqref{5} with the choice $g\equiv 1$, $\frac{\rr}{2}\le \tau_{1}<\tau_{2}\le \rr$, $\varphi\in C^{\infty}_{c}(B_{\rr})$ so that
\begin{flalign*}
\mathds{1}_{B_{\tau_{1}}}\le \varphi\le \mathds{1}_{B_{\tau_{2}}}\quad \mbox{and}\quad \snr{D\varphi}\le \frac{4}{\tau_{2}-\tau_{1}}
\end{flalign*}
and $\chi \in W^{1,\infty}(\mathbb{R},[0,1])$ with
\begin{flalign*}
\chi(t_{0}-\tau_{2}^{2})=0,\quad \chi\equiv 1 \ \ \mbox{on} \ \ (t_{0}-\tau_{1}^{2},t_{0}), \quad 0\le \partial_{t}\chi\le\frac{4}{(\tau_{2}-\tau_{1})^{2}}.
\end{flalign*}
Combining \eqref{5} with \eqref{unibd} we obtain
\begin{flalign}\label{6}
\sup_{t_{0}-\tau_{2}^{2}<t<t_{0}}\int_{B_{\rr}}&\varphi^{2}\chi\snr{Dv_{j}(x,t)}^{2} \ \dx+\int_{Q_{\rr}}\varphi^{2}\chi\snr{DV_{\ti{\mu},p}(Dv_{j})}^{2} \ \dx\nonumber \\
&+\varepsilon_{j}\int_{Q_{\rr}}\varphi^{2}\chi\snr{V_{\ti{\mu},2m\ti{q}}(Dv_{j})}^{2} \ \dy\nonumber \\
\le &\frac{c}{(\tau_{2}-\tau_{1})^{2}}\left(\int_{Q_{\tau_{2}}}\left[1+\ti{H}(Dv_{j})^{m\ti{q}} \ \dy\right]\right)^{\frac{1}{m}}+\frac{c\mathcal{C}_{f}}{(\tau_{2}-\tau_{1})^{2}},
\end{flalign}
with $c\equiv c(\texttt{data})$. Now we set
\begin{flalign}\label{tn}
\tilde{n}:=\begin{cases}
\ n\quad &\mbox{if} \ \ n>2\\
\ \mbox{any number in} \ \ \left(2,\min\left\{2\left(\frac{d}{d(q-p)+p}-1\right),\frac{2(d-2)}{d(q-p)+p}\right\}\right)\quad &\mbox{if} \ \ n=2 \ \ \mbox{and} \ \ \tilde{q}=q-\frac{p}{2}\\
\ \mbox{any number in} \ \ \left(2,\frac{2p(d-2)}{2d-pd+2p}\right)\quad &\mbox{if} \ \ n=2 \ \ \mbox{and} \ \ \tilde{q}=1
\end{cases}
\end{flalign}
and notice that, if $p\ge 2$ 
\begin{flalign*}
\ti{H}(z)^{\frac{p}{2}}\ge\snr{V_{\ti{\mu},p}(z)}^{2}\ge \snr{z}^{p}\quad \mbox{for all} \ \ z\in \mathbb{R}^{n},
\end{flalign*}
or, if $1<p<2$,
\begin{flalign*}
\ti{H}(z)^{\frac{p}{2}}\ge \snr{V_{\ti{\mu},p}(z)}^{2}\ge 2^{\frac{p-2}{2}}\snr{z}^{p}\quad \mbox{for all} \ \ z\in \mathbb{R}^{n} \ \ \mbox{with} \ \ \snr{z}\ge \ti{\mu}.
\end{flalign*}
On a fixed time slice we use H\"older inequality and \eqref{6} to bound
\begin{flalign*}
\int_{B_{\rr}}&\varphi^{2\left(1+\frac{2}{\tilde{n}}\right)}\snr{Dv_{j}}^{p+\frac{4}{\tilde{n}}} \ \dx\le \left(\int_{B_{\rr}}\varphi^{\frac{2\tilde{n}}{\tilde{n}-2}}\snr{Dv_{j}}^{\frac{\tilde{n}p}{\tilde{n}-2}} \ \dx\right)^{\frac{\tilde{n}-2}{\tilde{n}}}\left(\int_{B_{\rr}}\varphi^{2}\snr{Dv_{j}}^{2} \ \dx\right)^{\frac{2}{\tilde{n}}}\nonumber \\
\le &c\left[\left(\int_{B_{\rr}}\varphi^{\frac{2\tilde{n}}{\tilde{n}-2}} \ \dx\right)^{\frac{\tilde{n}-2}{\tilde{n}}}+\left(\int_{B_{\rr}}\varphi^{\frac{2\tilde{n}}{\tilde{n}-2}}\snr{V_{\ti{\mu},p}(Dv_{j})}^{\frac{2\tilde{n}}{\tilde{n}-2}} \ \dx\right)^{\frac{\tilde{n}-2}{\tilde{n}}}\right]\left(\int_{B_{\rr}}\varphi^{2}\snr{Dv_{j}}^{2} \ \dx\right)^{\frac{2}{\tilde{n}}}\nonumber \\
\le &c\left[\int_{B_{\rr}}\snr{D\varphi}^{2} \ \dx+\int_{B_{\rr}}\snr{D(\varphi V_{\ti{\mu},p}(Dv_{j}))}^{2} \ \dx\right]\left(\int_{B_{\rr}}\varphi^{2}\snr{Dv_{j}}^{2} \ \dx\right)^{\frac{2}{\tilde{n}}}\nonumber \\
\le &c\left[\int_{B_{\rr}}\snr{D\varphi}^{2} \ \dx+\int_{B_{\rr}}\varphi^{2}\snr{DV_{\ti{\mu},p}(Dv_{j})}^{2} \ \dx+\int_{B_{\rr}}\snr{V_{\ti{\mu},p}(Dv_{j})}^{2}\snr{D\varphi}^{2} \ \dx\right]\left(\int_{B_{\rr}}\varphi^{2}\snr{Dv_{j}}^{2} \ \dx\right)^{\frac{2}{\tilde{n}}}.
\end{flalign*}
We multiply both sides of the inequality in the previous display by $\chi^{1+\frac{2}{\tilde{n}}}$, integrate in time for $t\in (t_{0}-\tau_{2}^{2},t_{0})$, take the supremum in the time variable of the last integral on the right-hand side, use \eqref{6} and eventually get
\begin{flalign}\label{7}
\int_{Q_{\tau_{1}}}&\snr{Dv_{j}}^{p+\frac{4}{\tilde{n}}} \ \dy\le \frac{c}{(\tau_{2}-\tau_{1})^{2\left(1+\frac{2}{\tilde{n}}\right)}}\left(\int_{Q_{\tau_{2}}}\left[1+\ti{H}(Dv_{j})^{m\ti{q}}\right] \ \dy\right)^{\frac{1}{m}\left(1+\frac{2}{\tilde{n}}\right)}+\frac{c}{(\tau_{2}-\tau_{1})^{2\left(1+\frac{2}{\tilde{n}}\right)}}\nonumber\\
\le &\frac{c}{(\tau_{2}-\tau_{1})^{2\left(1+\frac{2}{\tilde{n}}\right)}}\left(\int_{Q_{\tau_{2}}}\snr{Dv_{j}}^{2m\ti{q}} \ \dy\right)^{\frac{1}{m}\left(1+\frac{2}{\tilde{n}}\right)}+\frac{c}{(\tau_{2}-\tau_{1})^{2\left(1+\frac{2}{\tilde{n}}\right)}},
\end{flalign}
where $c\equiv c(\texttt{data},\mathcal{C}_{f})$. We can rearrange \eqref{7} in the following way:
\begin{flalign}\label{8}
\nr{Dv_{j}}_{L^{p+\frac{4}{\tilde{n}}}(B_{\tau_{1}}\times (t_{0}-\tau_{1}^{2},t_{0}))}\le \frac{c}{(\tau_{2}-\tau_{1})^{\frac{2(\ti{n}+2)}{\ti{n}p+4}}}\nr{Dv_{j}}_{L^{2m\ti{q}}(B_{\tau_{2}}\times (t_{0}-\tau_{2}^{2},t_{0}))}^{\frac{2\ti{q}(\ti{n}+2)}{\ti{n}p+4}}+ \frac{c}{(\tau_{2}-\tau_{1})^{\frac{2(\ti{n}+2)}{\ti{n}p+4}}}.
\end{flalign}
Notice that, by \eqref{gamma} and \eqref{pq}, there holds that 
\begin{flalign}\label{pq.1}
p\le q<2m\ti{q}<p+\frac{4}{\ti{n}},
\end{flalign}
so we can apply the interpolation inequality 
\begin{flalign}\label{inter}
\nr{Dv_{j}}_{L^{2m\ti{q}}(B_{\tau_{2}}\times (t_{0}-\tau_{2}^{2},t_{0}))}\le \nr{Dv_{j}}_{L^{p}(B_{\tau_{2}}\times (t_{0}-\tau_{2}^{2},t_{0}))}^{1-\theta}\nr{Dv_{j}}_{L^{p+\frac{4}{\ti{n}}}(B_{\tau_{2}}\times (t_{0}-\tau_{2}^{2},t_{0}))}^{\theta},
\end{flalign}
where $\theta\in (0,1)$ solves
\begin{flalign*}
\frac{1}{2m\ti{q}}=\frac{1-\theta}{p}+\frac{\ti{n}\theta}{\ti{n}p+4}\ \ \Rightarrow \ \ \theta=\frac{(\ti{n}p+4)(2m\ti{q}-p)}{8m\ti{q}}.
\end{flalign*}
Plugging \eqref{inter} into \eqref{8} we get
\begin{flalign}\label{9}
\nr{Dv_{j}}_{L^{p+\frac{4}{\tilde{n}}}(B_{\tau_{1}}\times (t_{0}-\tau_{1}^{2},t_{0}))}\le &\frac{c}{(\tau_{2}-\tau_{1})^{\frac{2(\ti{n}+2)}{\ti{n}p+4}}}\nr{Dv_{j}}_{L^{p+\frac{4}{\tilde{n}}}(B_{\tau_{2}}\times(t_{0}-\tau_{2}^{2},t_{0}))}^{\frac{2\theta\ti{q}(\ti{n}+2)}{\ti{n}p+4}}\nr{Dv_{j}}_{L^{p}(B_{\tau_{2}}\times(t_{0}-\tau_{2}^{2},t_{0}))}^{\frac{2(1-\theta)\ti{q}(\ti{n}+2)}{\ti{n}p+4}}\nonumber \\
&+\frac{c}{(\tau_{2}-\tau_{1})^{\frac{2(\ti{n}+2)}{\ti{n}p+4}}},
\end{flalign}
with $c\equiv c(\texttt{data},\mathcal{C}_{f})$. By \eqref{pq} and \eqref{tn} there holds that
\begin{flalign*}
\frac{2\theta\ti{q}(\ti{n}+2)}{\ti{n}p+4}<1,
\end{flalign*}
so we can apply Young inequality with conjugate exponents
\begin{flalign}\label{conex}
\left(\frac{4m}{(2m\ti{q}-p)(\ti{n}+2)},\frac{4m}{4m-(2m\ti{q}-p)(\ti{n}+2)}\right)
\end{flalign}
to get
\begin{flalign}\label{lip4}
\nr{Dv_{j}}_{L^{p+\frac{4}{\tilde{n}}}(B_{\tau_{1}}\times (t_{0}-\tau_{1}^{2},t_{0}))}\le &\frac{1}{2}\nr{Dv_{j}}_{L^{p+\frac{4}{\tilde{n}}}(B_{\tau_{2}}\times (t_{0}-\tau_{2}^{2},t_{0}))}\nonumber \\
&+\frac{c(\texttt{data},\mathcal{C}_{f})}{(\tau_{2}-\tau_{1})^{\hat{\theta}}}\left[1+\nr{Dv_{j}}_{L^{p}(B_{\tau_{2}}\times (t_{0}-\tau_{2}^{2},t_{0}))}^{\beta}\right],
\end{flalign}
where we set $\hat{\theta}:=\frac{8m(\ti{n}+2)}{(\ti{n}p+4)[4m-(2m\ti{q}-p)(\ti{n}+2)]}$ and $\beta:=\frac{8m(1-\theta)\ti{q}(\ti{n}+2)}{[4m-(2m\ti{q}-p)(\ti{n}+2)](\ti{n}p+4)}$. Now we are in position to apply Lemma \ref{l0} and \eqref{unibd} to the inequality in the previous display and conclude with
\begin{flalign}
&\nr{Dv_{j}}_{L^{p+\frac{4}{\tilde{n}}}(B_{\rr/2}\times (t_{0}-(\rr/2)^{2},t_{0}))}\le \frac{c}{\rr^{\hat{\theta}}}\left[1+\nr{Dv_{j}}_{L^{p}(B_{\rr}\times (t_{0}-\rr^{2},t_{0}))}^{\beta}\right]\nonumber \\
&\qquad\le \frac{c}{\rr^{\hat{\theta}}}\left[\nr{Df}_{L^{r}(\Omega_{T})}^{r\beta}+\nr{\partial_{t}f}^{\beta p'}_{L^{p'}(0,T;W^{-1,p'}(\Omega))}+1\right]\label{10}
\end{flalign}
for $c\equiv c(\texttt{data})$. Finally, H\"older inequality and \eqref{10} in particular imply that
\begin{flalign}\label{11}
\nr{Dv_{j}}_{L^{s}(B_{\rr/2}\times (t_{0}-(\rr/2)^{2},t_{0}))}\le \frac{c(\texttt{data},\mathcal{C}_{f},s)}{\rr^{\hat{\theta}}}\quad \mbox{for all} \ \ s\in \left[1,p+\frac{4}{\ti{n}}\right],
\end{flalign}
thus \eqref{pq} and \eqref{tn} render that $s=q$ and $s=2m\ti{q}$ are both admissible choices. In the previous two displays, we also expanded the expression of $\mathcal{C}_{f}$. 
 \subsubsection*{Step 5: Fractional differentiability in space} Let $t_{0}\in (0,T)$ be any number and $\varphi\in C^{\infty}_{c}(B_{\rr})$ and $\chi\in W^{1,\infty}(\mathbb{R},[0,1])$ be two cut-off functions satisfying
 \begin{flalign}\label{phi}
 \mathds{1}_{B_{\rr/4}}\le \varphi\le \mathds{1}_{B_{\rr/2}}\quad \mbox{and}\quad \snr{D\varphi}\le \frac{4}{\rr}
 \end{flalign}
 and
 \begin{flalign}\label{chi}
 \chi(t_{0}-\rr^{2}/4)=0,\quad \chi=1 \ \ \mbox{on} \ \ (t_{0}-\rr^{2}/16,t_{0}),\quad 0\le \partial_{t}\chi\le \frac{4}{\rr^{2}}
 \end{flalign}
respectively. If $p\ge 2$, by Lemma \ref{l1} we have
\begin{flalign}\label{fr5}
\int_{Q_{\rr/2}}&\varphi^{2}\chi\snr{\Di V_{\ti{\mu},p}(Dv_{j})}^{2} \ \dy\sim \snr{h}^{-2}\int_{Q_{\rr/2}}\varphi^{2}\chi\mathcal{D}(h)^{\frac{p-2}{2}}\snr{\tau_{h}Dv_{j}}^{2} \ \dy\nonumber \\
\gtrsim& \snr{h}^{-2}\int_{Q_{\rr/2}}\varphi^{2}\chi\snr{\tau_{h}Dv_{j}}^{p} \ \dy,
\end{flalign}
while, for $1<p<2$ we have that
\begin{flalign}\label{fr6}
\snr{h}^{-p}\int_{Q_{\rr/2}}&\varphi^{2}\chi\snr{\tau_{j}Dv_{j}}^{p} \ \dy\le \left(  \snr{h}^{-2}\int_{Q_{\rr/2}}\varphi^{2}\chi\mathcal{D}(h)^{\frac{p-2}{2}}\snr{\tau_{h}Dv_{j}}^{2} \ \dx\right)^{\frac{p}{2}}\left(\int_{Q_{\rr/2}}\varphi^{2}\chi\mathcal{D}(j)^{\frac{p}{2}} \ \dy\right)^{\frac{2-p}{2}}\nonumber \\
\lesssim &\left(\snr{h}^{-2}\int_{Q_{\rr/2}}\varphi^{2}\chi\snr{\Di V_{\ti{\mu},p}(Dv_{j})}^{2} \ \dy\right)^{\frac{p}{2}}\left(\int_{Q_{\rr/2}}\varphi^{2}\chi\mathcal{D}(h)^{\frac{p}{2}} \ \dy\right)^{\frac{2-p}{2}}.
\end{flalign}
Therefore, if $p\ge 2$, by \eqref{fr5}, \eqref{2} with $g\equiv 1$, $\varphi$ and $\chi$ as in \eqref{phi}-\eqref{chi}, \eqref{6} and \eqref{11} we obtain
\begin{flalign}\label{fr7}
&\left(\limsup_{\snr{h}\to 0}\int_{Q_{\rr/4}}\left| \ \frac{\tau_{h}Dv_{j}}{\snr{h}^{\frac{2}{p}}} \ \right|^{p} \ \dy\right)\lesssim \limsup_{\snr{h}\to 0}\int_{Q_{\rr/2}}\varphi^{2}\chi\snr{\Di V_{\ti{\mu},p}(Dv_{j})}^{2} \ \dy\nonumber \\
&\qquad \lesssim \limsup_{\snr{h}\to 0}\mathcal{I}(h)\lesssim \rr^{-2}\left[1+\left(\int_{Q_{\rr/2}}\ti{H}(Dv_{j})^{m\ti q} \ \dy\right)^{\frac{1}{m}}\right]\lesssim \rr^{-\tilde{\theta}},
\end{flalign}
while, for $1<p<2$ we have, using also \eqref{unibd}
\begin{flalign}\label{fr8}
&\limsup_{\snr{h}\to 0}\left(\int_{Q_{\rr/4}}\left| \ \frac{\tau_{h}Dv_{j}}{\snr{h}} \ \right|^{p}\right)\lesssim \left(\limsup_{\snr{h}\to 0}\mathcal{I}(h)\right)^{\frac{p}{2}}\mathcal{C}_{f}^{\frac{2-p}{2}}\nonumber \\
&\qquad \lesssim\rr^{-p}\left[1+\left(\int_{Q_{\rr/2}}\ti{H}(Dv_{j})^{m\ti{q}} \ \dy\right)^{\frac{1}{m}}\right]^{\frac{p}{2}}\lesssim \rr^{-\tilde{\theta}},
\end{flalign}
In both, \eqref{fr7}-\eqref{fr8}, $\tilde{\theta}\equiv \tilde{\theta}(n,p,q,d)$ and the constants implicit in "$\lesssim$" depend on $(\texttt{data},\mathcal{C}_{f})$. Combining \eqref{fr7}-\eqref{fr8}, Proposition \ref{fracsob} and a standard covering argument, we can conclude that 
\begin{flalign}\label{fr9}
Dv_{j}\in L^{p}_{\loc}(0,T;W^{\varsigma,p}_{\loc}(\Omega,\mathbb{R}^{n}))\quad \mbox{for all} \ \ \varsigma\in \left(0,\min\left\{1,\frac{2}{p}\right\}\right).
\end{flalign}
\subsubsection*{Step 6: Fractional differentiability in time}
We aim to prove that 
\begin{flalign}\label{fr1}
\taa_{j}(\cdot,\cdot,Dv_{j})\in L^{l}_{\loc}(0,T;W^{1,l}_{\loc}(\Omega,\mathbb{R}^{n}))\quad \mbox{for some} \ \ l\equiv l(n,p,q,d)\in (1,\min\{2,p\}).
\end{flalign}
The forthcoming argument appears for instance in \cite{dumi} for the $p$-laplacean case with $p\ge 2$. Before going on, let us record some computations which will be helpful in a few lines. By the definition given in \eqref{r} it is clear that
\begin{flalign}\label{fr3.1}
\max\left\{\frac{p}{2},q-\frac{p}{2},m\ti{q}\right\}=m\ti{q}.
\end{flalign} 
Moreover, by \eqref{pq.1} we also have that there exists $l\in (1,2)$ so that
\begin{flalign}\label{fr2.1}
\max\left\{s_{1},s_{2}\right\}<p+\frac{4}{\ti{n}},
\end{flalign}
where we set 
\begin{flalign*}
s_{1}:=\frac{2l(m\ti{q}-1)}{2-l}\qquad\mbox{and}\qquad s_{2}:=\frac{dl(q-1)}{(d-l)}.
\end{flalign*}
For $\varphi,\chi$ as in \eqref{phi}-\eqref{chi} and $h$ as in \eqref{hhh}, we expand
\begin{flalign*}
\int_{Q_{\rr/2}}&\left[\varphi^{2}\chi\snr{\tau_{h} \taa_{j}(\cdot,t,Dv_{j})}\right]^{l} \ \dy\lesssim \int_{Q_{\rr/2}}\left[\varphi^{2}\chi\snr{\taa_{j}(x+h,t,Dv_{j}(x+h))-\taa_{j}(x,t,Dv_{j}(x+h))}\right]^{l} \ \dy\nonumber \\
&+\int_{Q_{\rr/2}}\left[\varphi^{2}\chi\snr{\taa_{j}(x,t,Dv_{j}(x+h))-\taa_{j}(x,t,Dv_{j}(x))}\right]^{l} \ \dy=:\mbox{(I)}+\mbox{(II)}
\end{flalign*}
and estimate, via $\eqref{refreg}_{3}$, \eqref{fr2.1} and H\"older inequality,
\begin{flalign*}
\mbox{(I)}\lesssim& \snr{h}^{l}\int_{Q_{\rr/2}}\varphi^{2l}\chi^{l}\left(\int_{0}^{1}\gamma(x+h\lambda,t) \ \d\lambda\right)^{l}\left[1+\mathcal{D}(h)^{\frac{l(q-1)}{2}}\right] \ \dy\nonumber \\
\lesssim&\snr{h}^{l}\nr{\gamma}_{L^{d}(\Omega_{T})}^{l}\left(\int_{Q_{\rr/2}}\left[1+\mathcal{D}(h)^{\frac{s_{2}}{2}}\right] \ \dy\right)^{\frac{l(q-1)}{s_{2}}}.
\end{flalign*}
Concerning term $\mbox{(II)}$ we distinguish three cases: $q\ge p\ge 2$, $q\ge 2>p$ and $2>q\ge p$. If $q\ge p\ge 2$, via $\eqref{refreg}_{1,3}$, \eqref{fr3.1}, \eqref{fr2.1}, H\"older inequality, Lemmas \ref{l1} and \ref{l6} we get
\begin{flalign*}
\mbox{(II)}\lesssim&\int_{Q_{\rr/2}}\varphi^{2l}\chi^{l}\left[\mathcal{D}(h)^{\frac{p-2}{2}}+\mathcal{D}(h)^{\frac{q-2}{2}}\right]^{l}\snr{\tau_{j}Dv_{j}}^{l} \ \dy\nonumber \\
&+\varepsilon_{j}\int_{Q_{\rr/2}}\left[\mathcal{D}(h)^{\frac{2m\ti{q}-2}{2}}\snr{\tau_{h}Dv_{j}}\right]^{l} \ \dy\nonumber \\
\lesssim&\snr{h}^{l}\left(\int_{Q_{\rr/2}}\varphi^{2l}\chi^{l}\mathcal{D}(h)^{\frac{l(p-2)}{2-l}} \ \dy\right)^{\frac{2-l}{2}}\left(\int_{Q_{\rr/2}}\varphi^{2l}\chi^{l}\snr{\Di V_{\ti{\mu},p}(Dv_{j})}^{2} \ \dy\right)^{\frac{l}{2}}\nonumber \\
&+\snr{h}^{l}\left(\int_{Q_{\rr/2}}\varphi^{2l}\chi^{l}\mathcal{D}(h)^{\frac{l(2q-p-2)}{2(2-l)}} \ \dy\right)^{\frac{2-l}{2}}\left(\int_{Q_{\rr/2}}\varphi^{2l}\chi^{l}\snr{\Di V_{\ti{\mu},p}(Dv_{j})}^{2} \ \dy\right)^{\frac{l}{2}}\nonumber \\
&+\snr{h}^{l}\left(\varepsilon_{j}\int_{Q_{\rr/2}}\varphi^{2l}\chi^{l}\mathcal{D}(h)^{\frac{l(m\ti{q}-1)}{2-l}} \ \dy\right)^{\frac{2-l}{2}}\left(\varepsilon_{j}\int_{Q_{\rr/2}}\varphi^{2l}\chi^{l}\snr{V_{\ti{\mu},2m\ti{q}}(Dv_{j})}^{2} \ \dy\right)^{\frac{l}{2}}.
\end{flalign*}
For $q\ge 2>p$, recalling \eqref{extra} we obtain 
\begin{flalign*}
\mbox{(II)}\lesssim& \snr{h}^{l}\mu^{p-2}\left(\int_{Q_{\rr/2}}\varphi^{2l}\chi^{l}\snr{\Di Dv_{j}}^{p} \ \dy\right)^{\frac{l}{p}}\nonumber \\
&+\snr{h}^{l}\left(\int_{Q_{\rr/2}}\varphi^{2l}\chi^{l}\mathcal{D}(h)^{\frac{l(2q-p-2)}{2(2-l)}} \ \dy\right)^{\frac{2-l}{2}}\left(\int_{Q_{\rr/2}}\varphi^{2l}\chi^{l}\snr{\Di V_{\ti{\mu},p}(Dv_{j})}^{2} \ \dy\right)^{\frac{l}{2}}\nonumber \\
&+\snr{h}^{l}\left(\varepsilon_{j}\int_{Q_{\rr/2}}\varphi^{2l}\chi^{l}\mathcal{D}(h)^{\frac{l(m\ti{q}-1)}{2-l}} \ \dy\right)^{\frac{2-l}{2}}\left(\varepsilon_{j}\int_{Q_{\rr/2}}\varphi^{2l}\chi^{l}\snr{V_{\ti{\mu},2m\ti{q}}(Dv_{j})}^{2} \ \dy\right)^{\frac{l}{2}}.
\end{flalign*}
Finally, when $2>q\ge p$ we use \eqref{extra} to conclude that
\begin{flalign*}
\mbox{(II)}\lesssim&\snr{h}^{l}\mu^{p-2}\left(\int_{Q_{\rr/2}}\varphi^{2l}\chi^{l}\snr{\Di Dv_{j}}^{p} \ \dy\right)^{\frac{l}{p}}\nonumber \\
&+\int_{Q_{\rr/2}\cap\{\mathcal{D}(h)\le 1\}}\varphi^{2l}\chi^{l}\left[\mathcal{D}(h)^{\frac{q-p}{2}}\mathcal{D}(h)^{\frac{p-2}{2}}\snr{\tau_{h}Dv_{j}}\right]^{l} \ \dy\nonumber \\
&+\int_{Q_{\rr/2}\cap\{\mathcal{D}(h)>1\}}\left[\mathcal{D}(h)^{\frac{q-2}{2}}\snr{\tau_{h}Dv_{j}}\right]^{l} \ \dy\nonumber \\
&+\snr{h}^{l}\left(\varepsilon_{j}\int_{Q_{\rr/2}}\varphi^{2l}\chi^{l}\mathcal{D}(h)^{\frac{l(m\ti{q}-1)}{2-l}} \ \dy\right)^{\frac{2-l}{2}}\left(\varepsilon_{j}\int_{Q_{\rr/2}}\varphi^{2l}\chi^{l}\snr{V_{\ti{\mu},2m\ti{q}}(Dv_{j})}^{2} \ \dy\right)^{\frac{l}{2}}\nonumber \\
\lesssim&\snr{h}^{l}(\mu^{p-2}+1)\left(\int_{Q_{\rr/2}}\varphi^{2l}\chi^{l}\snr{\Di Dv_{j}}^{p} \ \dy\right)^{\frac{l}{p}}\nonumber \\
&+\snr{h}^{l}\left(\varepsilon_{j}\int_{Q_{\rr/2}}\varphi^{2l}\chi^{l}\mathcal{D}(h)^{\frac{l(m\ti{q}-1)}{2-l}} \ \dy\right)^{\frac{2-l}{2}}\left(\varepsilon_{j}\int_{Q_{\rr/2}}\varphi^{2l}\chi^{l}\snr{V_{\ti{\mu},2m\ti{q}}(Dv_{j})}^{2} \ \dy\right)^{\frac{l}{2}}.
\end{flalign*}
Merging the content of all the previous displays and using Lemma \ref{diffquo}, \eqref{6} with $\tau_{1},\tau_{2}$ replaced by $\frac{\rr}{4},\frac{\rr}{2}$ respectively and \eqref{11}, we obtain
\begin{flalign*}
\limsup_{\snr{h}\to 0}&\int_{Q_{\rr/4}}\snr{\Di \taa_{j}(x,t,Dv_{j})}^{l}\ \dy\lesssim \nr{\gamma}_{L^{d}(\Omega_{T})}^{l}\left(\int_{Q_{\rr/2}}1+\snr{Dv_{j}}^{s_{2}} \ \dy\right)^{\frac{l(q-1)}{s_{2}}}\nonumber \\
&+(\ti{\mu}^{p-2}+1)\left(\limsup_{\snr{h}\to 0}\int_{Q_{\rr/2}}\snr{\Di Dv_{j}}^{p} \ \dy\right)^{\frac{l}{p}}\nonumber \\
&+\left(\int_{Q_{\rr/2}}\left[1+\snr{Dv_{j}}^{s_{1}}\right] \ \dy\right)^{\frac{2-l}{2}}\left(\limsup_{\snr{h}\to 0}\int_{Q_{\rr/2}}\snr{\Di V_{\ti{\mu},p}(Dv_{j})}^{2} \ \dy\right)^{\frac{l}{2}}\nonumber \\
&+\left(\varepsilon_{j}\int_{Q_{\rr/2}}\left[1+\snr{Dv_{j}}^{s_{1}}\right] \ \dy\right)^{\frac{2-l}{2}}\left(\limsup_{\snr{h}\to 0}\varepsilon_{j}\int_{Q_{\rr/2}}\snr{\Di V_{\ti{\mu},2m\ti{q}}(Dv_{j})}^{2} \ \dy\right)^{\frac{l}{2}}\lesssim \rr^{-\tilde{\theta}}.
\end{flalign*}
Finally, applying Fatou's lemma and Lemma \ref{diffquo} on the left-hand side of the chain of inequalities displayed above we obtain that
\begin{flalign}\label{fr9.1}
\int_{Q_{\rr/4}}\snr{D\taa_{j}(x,t,Dv_{j})}^{l} \ \dy\le c\rr^{-\tilde{\theta}},
\end{flalign}
with $c\equiv c(\texttt{data},\mathcal{C}_{f},\ti{\mu})$ and $\tilde{\theta}\equiv \tilde{\theta}(n,p,q,d)$. With \eqref{fr9.1} and a standard covering argument we deduce \eqref{fr1}. Now, whenever we consider a subset of type $\tilde{\Omega}\times(t_{1},t_{2})\Subset \Omega_{T}$ with $\tilde{\Omega}\Subset \Omega$ open, from \eqref{fr9.1} and \eqref{fr1} and a covering argument we have that
\begin{flalign}\label{fr11}
\nr{\diver \ \taa_{j}(\cdot,\cdot,Dv_{j})}_{L^{l}(\tilde{\Omega}\times (t_{1},t_{2}))}\le c\nr{D \taa_{j}(\cdot,\cdot,Dv_{j})}_{L^{l}(\tilde{\Omega}\times (t_{1},t_{2}))}\le c,
\end{flalign}
for $c\equiv c(\texttt{data},\mathcal{C}_{f},\mu,t_{1},T-t_{2},\dist(\tilde{\Omega},\partial \Omega))$. Finally, integrating by parts in \eqref{wfj} and using \eqref{fr11} we obtain that
\begin{flalign}\label{fr12}
\partial_{t}v_{j}\in L^{l}_{\loc}(\Omega_{T})\qquad \mbox{with} \ \ l\equiv l(n,p,q,d)\in (1,\min\{p,2\}).
\end{flalign}
\subsubsection*{Step 6: Convergence}  A standard covering argument combined with Proposition \ref{fracsob}, \eqref{11}, \eqref{fr7}-\eqref{fr9} and \eqref{fr11}-\eqref{fr12} respectively then implies that if $\tilde{\Omega}\Subset \Omega$ is any open subset and $(t_{1},t_{2})\Subset (0,T)$ is an interval, then
\begin{flalign}
&\nr{Dv_{j}}_{L^{s}(\tilde{\Omega}\times (t_{1},t_{2}))}\le c\quad \mbox{for all} \ \ s\in\left[1,p+\frac{4}{\ti{n}}\right];\label{13.1}\\
&\nr{v_{j}}_{L^{p}(t_{1},t_{2};W^{1+\varsigma}(\tilde{\Omega}))}\le c\quad \mbox{for all} \ \ \varsigma\in\left(0,\min\left\{1,\frac{2}{p}\right\}\right);\label{13}\\
&\nr{v_{j}}_{W^{\iota,l}(t_{1},t_{2};L^{l}(\tilde{\Omega}))}\le c\quad \mbox{for all} \ \ \iota\in (0,1),\label{13.2}
\end{flalign}
with $c\equiv c(\texttt{data},s,\varsigma,\iota,\mathcal{C}_{f},t_{1},T-t_{2},\dist(\tilde{\Omega},\partial \Omega))$. Estimates \eqref{13} and \eqref{13.2} render that
\begin{flalign*}
\{v_{j}\} \ \mbox{is uniformly bounded in} \ W^{\iota,l}_{\loc}(0,T;L^{l}_{\loc}(\Omega))\cap L^{p}_{\loc}(0,T;W^{1+\varsigma,p}_{\loc}(\Omega)),
\end{flalign*}
therefore we can first choose $\iota\in \left(\frac{p-l}{lp},1\right)$ so that $l>\frac{p}{1+\iota p}$ and then apply Lemma \ref{al} with $a_{1}=l$, $a_{2}=p$, $\sigma=\iota$, $X=W^{1+\varsigma,p}_{\loc}(\Omega)$, $B=W^{1,l}_{\loc}(\Omega)$ and $Y=L^{l}_{\loc}(\Omega)$ to conclude that 
\begin{flalign}\label{cv2}
\mbox{there exists a subsequence}\ \{v_{j}\} \ \mbox{strongly converging to} \ v \ \mbox{in} \ L^{l}_{\loc}(0,T;W^{1,l}_{\loc}(\Omega)),
\end{flalign}
where we also used that $l<p$. Using \eqref{13.1} we also see that, again up to subsequences,
\begin{flalign}\label{cv3}
Dv_{j}\rightharpoonup Dv\quad \mbox{in} \ \ L^{s}_{\loc}(\Omega_{T},\mathbb{R}^{n}) \quad \mbox{for all} \ \ s\in\left[1,p+\frac{4}{\ti{n}}\right]
\end{flalign}
which assures that
\begin{flalign}\label{cv5}
\nr{Dv}_{L^{s}(\tilde{\Omega}\times (t_{1},t_{2}))}\le c(\texttt{data},s,\mathcal{C}_{f},t_{1},T-t_{2},\dist(\tilde{\Omega},\partial \Omega))\quad \mbox{and}\quad \left.v\right|_{\partial_{par}\Omega_{T}}=\left.f\right|_{\partial_{par}\Omega_{T}}.
\end{flalign}
By \eqref{pq.1}, \eqref{cv2}, \eqref{cv3}, \eqref{cv5} and the interpolation inequality 
\begin{flalign*}
\nr{Dv_{j}-Dv}_{L^{s}(\tilde{\Omega}\times (t_{1},t_{2}))}\le& \nr{Dv_{j}-Dv}_{L^{l}(\tilde{\Omega}\times (t_{1},t_{2})))}^{\theta}\nr{Dv_{j}-Dv}_{L^{p+\frac{4}{\tilde{n}}}(\tilde{\Omega}\times (t_{1},t_{2})))}^{1-\theta}\nonumber \\
\le &c\nr{Dv_{j}-Dv}_{L^{l}(\tilde{\Omega}\times (t_{1},t_{2})))}^{\theta}
\end{flalign*}
with $c\equiv c(\texttt{data},s,\mathcal{C}_{f},t_{1},T-t_{2},\dist(\tilde{\Omega},\partial \Omega))$ and
\begin{flalign*}
\frac{1}{s}=\frac{\tilde{n}\theta}{\tilde{n}p+4}+\frac{1-\theta}{l} \ \Longrightarrow \ \theta=\frac{(\tilde{n}p+4)(s-l)}{s(\tilde{n}p+4-\tilde{n}l)},
\end{flalign*}
 we can conclude that
\begin{flalign}\label{cv4}
Dv_{j}\to Dv\quad \mbox{in} \ \ L^{s}_{\loc}(0,T;L^{s}_{\loc}(\Omega,\mathbb{R}^{n}))\quad \mbox{for all} \ \ s\in \left[1,p+\frac{4}{\tilde{n}}\right).
\end{flalign}
Once \eqref{cv4} is available, we can look back at \eqref{2} with $g\equiv 1$, send first $j\to \infty$ and then $\snr{h}\to 0$ and rearrange the right-hand side with the help of \eqref{10} to obtain \eqref{difff}. Moreover, using \eqref{cv4}, $\eqref{refreg}_{1}$ and the dominated convergence theorem, we can pass to the limit in \eqref{wfj} to conclude that $v$ satisfies
\begin{flalign}\label{cv6}
\int_{\Omega_{T}}\left[v\partial_{t}\varphi-a(x,t,Dv)\cdot D\varphi\right] \ \dy=0\quad \mbox{for all} \ \ \varphi\in C^{\infty}_{c}(\Omega_{T}).
\end{flalign}
Once \eqref{sss}, \eqref{cv6} and $\eqref{difff}$ are available, we can repeat the same computations leading to \eqref{fr1}-\eqref{fr12} with $\ti{a}(\cdot), v$ replacing $\taa_{j}(\cdot), v_{j}$ to obtain \eqref{cv7}.
\subsubsection*{Step 8: The initial boundary condition} 
With \eqref{cv6}, the energy estimate \eqref{unibd} and the continuity of $f$ in time prescribed by $\eqref{ggg}_{1}$, we can proceed exactly as in \cite[Section 6.5]{bdm} to verify the requirements of Definition \ref{d.1} (formulated for $v$ and $\ti{a}(\cdot)$ of course).
\section{Gradient bounds}\label{mose} This section is divided into two parts: in the first one we construct a sequence of maps satisfying suitable uniform estimates and in the second we prove that such a sequence converges to a weak solution of problem \eqref{pdd}.
\subsection{Uniform $L^{\infty}$-estimates}\label{inf}
We consider again Cauchy-Dirichlet problem \eqref{pdd} with $a(\cdot)$ described by \eqref{regg}-\eqref{pq} and $f$ as in \eqref{ggg}. To construct a suitable family of approximating problems, this time we only regularize the vector field $a(\cdot)$ in the gradient variable by convolution against a sequence $\{\phi_{j}\}$ of mollifiers of $\mathbb{R}^{n}$ with the following features:
\begin{flalign*}
\phi\in C^{\infty}_{c}(B_{1}), \quad \nr{\phi}_{L^{1}(\mathbb{R}^{n})}=1,\quad \phi_{j}(x):=j^{n}\phi(jx),\quad B_{3/4}\subset \supp (\phi).
\end{flalign*}
This leads to the definition of the approximating vector field
\begin{flalign}\label{lip3}
a_{j}(x,t,z):=\mint_{B_{1}}a(x,t,z+j^{-1}z')\phi(z') \ \dzz,
\end{flalign}
satisfying the structural conditions 
\begin{flalign}\label{reggej}
\begin{cases}
\ t\mapsto a_{j}(x,t,z)\quad &\mbox{measurable for all} \ \ x\in \Omega, z\in \mathbb{R}^{n}\\
\ x\mapsto a_{j}(x,t,z)\quad &\mbox{differentiable for all} \ \ t\in (0,T),z\in \mathbb{R}^{n}\\
\ z\mapsto a_{j}(x,t,z)\in C^{1}(\mathbb{R}^{n},\mathbb{R}^{n})\quad &\mbox{for all} \ \ (x,t)\in \Omega_{T}
\end{cases}
\end{flalign}
and
\begin{flalign}\label{refregj}
\begin{cases}
\ \snr{a_{j}(x,t,z)}+H_{j}(z)^{\frac{1}{2}}\snr{\partial_{z}a_{j}(x,t,z)}\le c\left[H_{j}(z)^{\frac{p-1}{2}}+H_{j}(z)^{\frac{q-1}{2}}\right]\\
\ \partial_{z}a_{j}(x,t,z)\ge cH_{j}(z)^{\frac{p-2}{2}}\snr{\xi}^{2}\\
\ \snr{\partial_{x}a_{j}(x,t,z)}\le c\gamma(x,t)\left[H_{j}(z)^{\frac{p-1}{2}}+H_{j}(z)^{\frac{q-1}{2}}\right],
\end{cases}
\end{flalign}
for all $(x,t)\in \Omega_{T}$, $z,\xi\in \mathbb{R}^{n}$, $\gamma$ as in \eqref{gamma}, with $c\equiv c(n,\nu,L,p,q)$, see \cite[Section 4.5]{dm} for more details on this matter. In \eqref{refregj}, 
\begin{flalign*}
\mu_{j}:=\mu+j^{-1}>0\quad \mbox{and}\quad H_{j}(z):=(\mu_{j}^{2}+\snr{z}^{2}).
\end{flalign*}
We then define problem
\begin{flalign}\label{pdjj}
\begin{cases}
\ \partial_{t}v-\diver \ a_{j}(x,t,Dv)=0\quad &\mbox{in} \ \ \Omega_{T}\\
\ v=f\quad &\mbox{on} \ \ \partial_{par}\Omega_{T},
\end{cases}
\end{flalign}
with $f$ as in \eqref{ggg}. By \eqref{reggej}-\eqref{refregj}, we see that the assumptions of Proposition \ref{p1} are satisfied, thus problem \eqref{pdjj} admits a solution $u_{j}\in L^{p}(0,T;W^{1,p}(\Omega))$ in the sense of Definition \ref{d.1}, satisfying \eqref{sss}, \eqref{cv7} and \eqref{difff}. In particular, \eqref{sss} authorizes to test \eqref{wfp} against test functions defined as products of $u_{j}$ with suitable cut-off functions, therefore, for such a solution, we can repeat almost the same computations leading to \eqref{5} (with $\varepsilon_{j}\equiv 0$, of course), for getting
\begin{flalign}\label{lip0}
&\frac{1}{2}\int_{B_{\rr}}\varphi^{2}\chi\left(\int_{0}^{\snr{D u_{j}}^{2}}g(s) \ \ds\right) \ \dx\nonumber \\
&\quad+\int_{Q_{\tau}}\varphi^{2}\chi\hat{\mathcal{G}}(Du_{j})\snr{D V_{\mu_{j},p}(Du_{j})}^{2} \ \dy\nonumber \\
&\quad\le c\left(\int_{Q_{\tau}}\chi^{m}\left(\snr{D\varphi}^{2m}+\varphi^{2m}\right)\hat{\mathcal{G}}(Du_{j})^{m}\left[1+H_{j}(Du_{j})^{m\left(q-\frac{p}{2}\right)}\right] \ \dy\right)^{\frac{1}{m}}\nonumber \\
&\quad+c\int_{Q_{\tau}}\left(\int_{0}^{\snr{Du_{j}}^{2}}g(s) \ \ds\right) \ \varphi^{2}\partial_{t}\chi \ \dy,
\end{flalign}
with $c\equiv c(\texttt{data})$, $g\in W^{1,\infty}(\mathbb{R})$ non-negative with bounded, non-negative, piecewise continuous first derivative, $\varphi\in C^{\infty}_{c}(B_{\rr},[0,1])$ and $\chi\in W^{1,\infty}([0,T])$. The quantity $\hat{\mathcal{G}}(Du_{j})$ is defined as in \emph{Step 3} of the proof of Proposition \ref{p1}, clearly with $u_{j}$ replacing $v_{j}$. For $i\in \N$, we inductively define radii $\rr_{i}:=\tau_{1}+(\tau_{2}-\tau_{1})2^{-i+1}$ with $\frac{\rr}{2}\le \tau_{1}<\tau_{2}\le \rr$,
select cut-off functions $\varphi_{i}\in C^{1}_{c}(B_{\rr})$ so that
\begin{flalign*}
&\mathds{1}_{B_{\rr_{i+1}}}\le \varphi_{i}\le \mathds{1}_{B_{\rr_{i}}}\quad \mbox{and}\quad \snr{D\varphi_{i}}\le \frac{4}{\rr_{i}-\rr_{i+1}}=\frac{2^{i+2}}{(\tau_{2}-\tau_{1})}
\end{flalign*}
and $\chi_{i}\in W^{1,\infty}_{0}((t_{0}-\rr^{2},t_{0}),[0,1])$ satisfying
\begin{flalign*}
&\chi_{i}(t_{0}-\rr_{j}^{2})=0,\quad \chi_{i}\equiv 1 \ \ \mbox{on} \ \ (t_{0}-\rr_{i+1}^{2},t_{0}),\quad  \snr{\partial_{t}\chi_{i}}\le \frac{4}{(\rr_{i}-\rr_{i+1})^{2}}\le \frac{2^{2i}}{(\tau_{2}-\tau_{1})^{2}}
\end{flalign*} 
and numbers
\begin{flalign}\label{ki}
\kk_{1}\equiv 0,\qquad \kk_{i}:=\frac{\Gamma}{m}+\omega \kk_{i-1} \ \ \mbox{for} \ \ i\ge 2,\qquad \alpha_{i}:=m\ti{q}+m\kk_{i},
\end{flalign}
where we set 
\begin{flalign}
\omega:=\frac{1}{m}\left[1+\frac{2}{\ti{n}}\right]\stackrel{\eqref{gamma}}{>}1\quad \mbox{and}\quad \Gamma:=\frac{p}{2}+\frac{2}{\tilde{n}}-m\ti{q}\stackrel{\eqref{pq}}{>}0.\label{s1}
\end{flalign}
In \eqref{lip0} we take $\varphi\equiv \varphi_{i}$, $\chi\equiv \chi_{i}$ and, for $M>0$ set
\begin{flalign*}
g(s)\equiv g_{i,M}(s):=\begin{cases}
\ (\mu_{j}^{2}+s)^{\kk_{i}}&\quad \mbox{if} \ \ s\le M\\
\ (\mu_{j}^{2}+M)^{\kk_{i}}&\quad \mbox{if} \ \ s> M,
\end{cases}
\end{flalign*}
which is admissible by construction in \eqref{lip0}. Clearly,
\begin{flalign}\label{lip6}
g_{i,M}(s)\le (\mu_{j}^{2}+s)^{\kk_{i}}\quad \mbox{for all} \ \ s\in [0,\infty).
\end{flalign}
All in all, \eqref{lip0} becomes
\begin{flalign}\label{lip1}
&\frac{1}{2}\int_{B_{\rr}}\varphi_{i}^{2}\chi_{i}\left(\int_{0}^{\snr{D u_{j}}^{2}}g_{i,M}(s) \ \ds\right) \ \dx\nonumber \\
&\quad+\int_{Q_{\tau}}\varphi_{i}^{2}\chi_{i}\hat{\mathcal{G}}_{i,M}(Du_{j})\snr{D V_{\mu_{j},p}(Du_{j})}^{2} \ \dy\nonumber \\
&\quad\le c\left(\int_{Q_{\tau}}\chi_{i}^{m}\left(\snr{D\ti{\varphi}}^{2m}+\varphi_{i}^{2m}\right)\hat{\mathcal{G}}_{i,M}(Du_{j})^{m}\left[1+H_{j}(Du_{j})^{m\left(q-\frac{p}{2}\right)}\right] \ \dy\right)^{\frac{1}{m}}\nonumber \\
&\quad+c\int_{Q_{\tau}}\left(\int_{0}^{\snr{Du_{j}}^{2}}g_{i,M}(s) \ \ds\right) \ \varphi_{i}^{2}\partial_{t}\chi_{i} \ \dy,
\end{flalign}
where we defined $\hat{\mathcal{G}}_{i,M}$ in the obvious way: $\hat{\mathcal{G}}_{i,M}(z):=\left(g_{i,M}(\snr{z}^{2})+\snr{z}^{2}g_{i,M}'(\snr{z}^{2})\right)$. As we only know that $\{u_{j}\}$ satisfies \eqref{sss}-\eqref{difff}, we have to proceed inductively. We shall prove that
\begin{flalign}\label{lip2}
H_{j}(Du_{j})^{\alpha_{i}}\in L^{1}(Q_{\rr_{i}})\Rightarrow H_{j}(Du_{j})^{\alpha_{i+1}}\in L^{1}(Q_{\rr_{i+1}})\quad \mbox{for all} \ \ i\in \mathbb{N}.
\end{flalign}
\subsubsection*{Basic step} Let us verify \eqref{lip2} for $i=1$. In this case, we immediately see that $\hat{\mathcal{G}}_{1,m}(Du_{j})\equiv 1$ and notice that, since the approximating sequence $\{u_{j}\}$ we choose satisfies \eqref{sss}-\eqref{difff}, all the computations made in \emph{Step 3} of Section \ref{high} are legal without further corrections to the growth of the vector field defined in \eqref{lip3}. Moreover, a quick inspection of estimates \eqref{5}-\eqref{6} points out the dependency of the constants from $\mathcal{C}_{f}$ is due only to the presence of the term multiplying $\varepsilon_{j}$, which, in the present case is zero. Hence, \eqref{lip1} becomes \eqref{6} with $\varepsilon_{j}\equiv 0$, $\varphi\equiv \varphi_{1}$ and $\chi\equiv \chi_{1}$. Since $\alpha_{1}=m\ti{q}$ and $\alpha_{2}=\frac{p}{2}+\frac{2}{\ti{n}}$, we can easily deduce from \eqref{7} (with $\tau_{1}=\rr_{2}$, $\tau_{2}=\rr_{1}$ and no dependencies of the constants from $\mathcal{C}_{f}$) that $H_{j}(Du_{j})^{\alpha_{2}}\in L^{1}(Q_{\rr_{2}})$. 
\subsubsection*{Induction step} We assume now that
\begin{flalign}\label{lip5}
H_{j}(Du_{j})^{\alpha_{i}}\in L^{1}(Q_{\rr_{i}})
\end{flalign}
and expand into \eqref{lip1} the expression of $\hat{\mathcal{G}}_{i,M}(Du_{j})$ for getting, after a few standard manipulations:
\begin{flalign*}
\frac{1}{2}\int_{B_{\rr_{i}}}&\varphi_{i}\chi_{i}\left(\int_{0}^{\min\{\snr{Du_{j}}^{2},M\}}(\mu_{j}^{2}+s)^{\kk_{i}} \ \ds\right) \ \dx\nonumber \\
&+\int_{Q_{\rr_{i}}\cap\{\snr{Du_{j}}^{2}
\le M\}}\varphi_{i}^{2}\chi_{i}(\mu_{j}^{2}+\snr{Du_{j}}^{2})^{\kk_{i}}\snr{DV_{\mu_{j},p}(Du_{j})}^{2} \ \dy\nonumber \\
\le &c(1+\kk_{i})\left(\int_{Q_{\rr_{i}}}\ti{\chi}^{m}\left(\snr{D\varphi_{i}}^{2m}+\varphi_{i}^{2m}\right)\left[1+H_{j}(Du_{j})^{m\left(\kk_{i}+q-\frac{p}{2}\right)}\right] \ \dy\right)^{\frac{1}{m}}\nonumber \\
&+\frac{c}{1+\kappa_{i}}\int_{Q_{\rr_{i}}}\varphi_{i}^{2}\partial_{t}\chi_{i}H_{j}(Du_{j})^{1+\kk_{i}} \ \dy\nonumber \\
\le &c(1+\kappa_{i})\left(\int_{Q_{\rr_{i}}}\left[\ti{\chi}^{m}\left(\snr{D\varphi_{i}}^{2m}+\varphi_{i}^{2m}\right)+\varphi_{i}^{2m}\snr{\partial_{t}\chi_{i}}^{m}\right]\left[1+H_{j}(Du_{j})^{m(\kk_{i}+\ti{q})}\right] \ \dy\right)^{\frac{1}{m}},
\end{flalign*}
for $c\equiv c(\texttt{data})$. For the inequality in the previous display we used in particular \eqref{lip6} and the definition of $\varphi_{i},\chi_{i}$. Now we can send $M\to \infty$ in the previous display and apply Fatou Lemma on the left-hand side, the dominated convergence theorem, $\eqref{ki}_{3}$ and \eqref{lip5} on the right-hand side to conclude with
\begin{flalign}\label{lip7}
\frac{1}{2}\int_{B_{\rr_{i}}}&\varphi_{i}\chi_{i}H_{j}(Du_{j})^{1+\kappa_{i}} \ \dx+(1+\kappa_{i})\int_{Q_{\rr_{i}}}\varphi_{i}^{2}\chi_{i}H_{j}(Du_{j})^{\kk_{i}}\snr{DV_{\mu_{j},p}(Du_{j})}^{2} \ \dy\nonumber \\
\le &c(1+\kappa_{i})^{2}\left(\int_{Q_{\rr_{i}}}\left[\chi_{i}^{m}\left(\snr{D\varphi_{i}}^{2m}+\varphi_{i}^{2m}\right)+\varphi_{i}^{2m}\snr{\partial_{t}\chi_{i}}^{m}\right]\left[1+H_{j}(Du_{j})^{\alpha_{i}}\right] \ \dy\right)^{\frac{1}{m}},
\end{flalign}
where $c\equiv c(\texttt{data})$. Next, with \eqref{difff} at hand, we compute
\begin{flalign*}
\snr{DH_{j}(Du_{j})^{\frac{p+2\kk_{i}}{4}}}^{2}=\left(\frac{p+2\kk_{i}}{p}\right)^{2}H_{j}(Du_{j})^{\kk_{i}}\snr{DH_{j}(Du_{j})^{\frac{p}{4}}}^{2}
\end{flalign*}
and 
\begin{flalign*}%\label{60}
\snr{DV_{\mu_{j},p}(Du)}^{2}=&\left(\frac{p-2}{2}\right)^{2}H_{j}(Du_{j})^{\frac{p-6}{2}}\snr{Du_{j}\cdot D^{2}u_{j}}^{2}\snr{Du_{j}}^{2}\nonumber \\
&+H_{j}(Du_{j})^{\frac{p-2}{2}}\snr{D^{2}u_{j}}^{2}+(p-2)H_{j}(Du_{j})^{\frac{p-4}{2}}\snr{Du_{j}\cdot D^{2}u_{j}}^{2}\nonumber \\
\ge &\min\{1,(p-1)\}H_{j}(Du_{j})^{\frac{p-2}{2}}\snr{D^{2}u_{j}}^{2},
\end{flalign*}
so, keeping in mind that 
\begin{flalign*}
\snr{DH_{j}(Du_{j})^{\frac{p}{4}}}^{2}\le \left(\frac{p}{2}\right)^{2}H_{j}(Du_{j})^{\frac{p-2}{2}}\snr{D^{2}u_{j}}^{2}
\end{flalign*}
we end up with
\begin{flalign}\label{19}
\snr{DH_{j}(Du_{j})^{\frac{p+2\kk_{i}}{4}}}^{2}\le  \frac{(p+2\kk_{i})^{2}}{4\min\{p-1,1\}}H_{j}(Du_{j})^{\kk_{i}}\snr{DV_{\mu_{j},p}(Du_{j})}^{2}.
\end{flalign}
Plugging \eqref{19} into \eqref{lip7} we obtain, after routine calculation,
\begin{flalign}\label{lip8}
&\sup_{t_{0}-(r_{i,2})^{2}<t<t_{0}}\int_{B_{\rr_{i}}}\varphi_{i}^{2}\chi_{i}H_{j}(Du_{j})^{1+\kappa_{i}} \ \dx+\int_{Q_{\rr_{i}}}\chi_{i}\snr{D(\varphi_{i}[H_{j}(Du_{j})^{\frac{p+2\kk_{i}}{4}}+1])}^{2} \ \dy\nonumber \\
&\qquad\le \sup_{t_{0}-(r_{i,2})^{2}<t<t_{0}}\int_{B_{\rr_{i}}}\varphi_{i}\chi_{i}H_{j}(Du_{j})^{1+\kappa_{i}} \ \dx\nonumber \\
&\qquad +c\int_{Q_{\rr_{i}}}\chi_{i}\left[\varphi_{i}^{2}\snr{DH_{j}(Du_{j})^{\frac{p+2\kk_{i}}{4}}}^{2}+\snr{D\ti{\varphi}}^{2}\left(H_{j}(Du_{j})^{\frac{p+2\kk_{i}}{2}}+1\right)\right] \ \dy\nonumber \\
&\qquad\le c(1+\kappa_{i})^{4}\left(\int_{Q_{\rr_{i}}}\left[\chi^{m}\left(\snr{D\varphi_{i}}^{2m}+\varphi_{i}^{2m}\right)+\varphi_{i}^{2m}\snr{\partial_{t}\chi_{i}}^{m}\right]\left[1+H_{j}(Du_{j})^{\alpha_{i}}\right] \ \dy\right)^{\frac{1}{m}},
\end{flalign}
with $c\equiv c(\texttt{data})$. For $\ti{n}$ as in \eqref{tn}, we define $\ti{\sigma}_{i}:=2(1+\kk_{i})\ti{n}^{-1}$. On a fixed time slice, we apply in sequence H\"older and Sobolev-Poincar\'e inequalities to get
\begin{flalign}\label{lip9}
\int_{B_{\rr_{i}}}&\varphi_{i}^{2\left(1+\frac{2}{\tilde{n}}\right)}H_{j}(Du_{j})^{\frac{p+2\kk_{i}}{2}+\ti{\sigma}_{i}} \ \dx\nonumber \\
\le& \left(\int_{B_{\rr_{i}}}\left[\varphi_{i}^{2}(H_{j}(Du_{j})^{\frac{p+2\kk_{i}}{2}}+1)\right]^{\frac{\tilde{n}}{\tilde{n}-2}} \ \dx\right)^{\frac{\tilde{n}-2}{\tilde{n}}}\left(\int_{B_{\rr_{i}}}\varphi_{i}^{2}H_{j}(Du_{j})^{\ti{\sigma}_{i}\frac{\tilde{n}}{2}} \ \dx\right)^{\frac{2}{\tilde{n}}}\nonumber \\
\le &c\left(\int_{B_{\rr_{i}}}\snr{D[\varphi_{i}(H_{j}(Du_{j})^{\frac{p+2\kk_{i}}{4}}+1)]}^{2} \ \dx\right)\left(\int_{B_{\rr_{i}}}\varphi_{i}^{2}H_{j}(Du_{j})^{\ti{\sigma}_{i}\frac{\tilde{n}}{2}} \ \dx\right)^{\frac{2}{\tilde{n}}},
\end{flalign}
for $c\equiv c(n,p,q,d)$. Now we multiply both sides of \eqref{lip9} by $\chi_{i}^{\frac{\tilde{n}+2}{\tilde{n}}}$, integrate with respect to $t\in (t_{0}-(r_{i,2})^{2},t_{0})$, take the supremum over $t\in (t_{0}-(r_{i,2})^{2},t_{0})$ on the right-hand side, use \eqref{lip8} and eventually obtain
\begin{flalign}\label{lip10}
\int_{Q_{\rr_{i}}}&(\varphi_{i}^{2}\chi_{i})^{1+\frac{2}{\tilde{n}}}\left[1+H_{j}(Du_{j})^{\frac{p+2\kk_{i}}{2}+\ti{\sigma}_{i}}\right] \ \dy\nonumber \\
\le &c(1+\kk_{i})^{4\left(1+\frac{2}{\tilde{n}}\right)}\left(\int_{Q_{\rr_{i}}}\left[\chi_{i}^{m}\left(\varphi_{i}^{2m}+\snr{D\varphi_{i}}^{2m}\right)+\varphi_{i}^{2m}\snr{\partial_{t}\chi_{i}}^{m}\right]\left[1+H_{j}(Du_{j})^{\alpha_{i}}\right] \ \dy\right)^{\frac{1}{m}\left(1+\frac{2}{\tilde{n}}\right)},
\end{flalign}
where $c\equiv c(\texttt{data})$. In the light of \eqref{ki}-\eqref{s1} we have
\begin{flalign}\label{lip11}
\frac{p}{2}+\kk_{i}+\ti{\sigma}_{i}=&\frac{p}{2}+\frac{2}{\ti{n}}+m\omega\kk_{i}=\left(\frac{p}{2}+\frac{2}{\ti{n}}-m\ti{q}\right)+m\left(\ti{q}+\omega\kk_{i}\right)\nonumber \\
=&m\left(\frac{\Gamma}{m}+\ti{q}+\omega\kk_{i}\right)=m\left(\ti{q}+\kk_{i+1}\right)=\alpha_{i+1},
\end{flalign}
so, recalling also the definition of $\chi_{i},\varphi_{i}$ \eqref{lip10} becomes
\begin{flalign*}
\int_{Q_{\rr_{i+1}}}H_{j}(Du_{j})^{\alpha_{i+1}} \ \dy\le \frac{c(\texttt{data},i)}{(\rr_{i}-\rr_{i+1})^{2}}\left(\int_{Q_{\rr_{i}}}\left[1+H_{j}(Du_{j})^{\alpha_{i}}\right] \ \dy\right)^{\frac{1}{m}\left(1+\frac{2}{\tilde{n}}\right)}\stackrel{\eqref{lip5}}{<}\infty
\end{flalign*}
and \eqref{lip5} is proved for all $i\in \mathbb{N}$.\\\\
Now we know that the quantity appearing on the right-hand side of \eqref{lip10} is finite for all $i\in \mathbb{N}$, we define
\begin{flalign*}
A_{i}:=\left(\mint_{Q_{\rr_{j}}}\left[1+H_{j}(Du_{j})^{\alpha_{i}}\right] \ \d z\right)^{\frac{1}{\alpha_{i}}}.
\end{flalign*}
From the definitions in \eqref{ki}, it is easy to see that whenever $i\ge 2$
\begin{flalign*}
\kk_{i}=\frac{\Gamma}{m}\sum_{l=0}^{i-2}\omega^{i}\quad \mbox{and}\quad \alpha_{i}=m\ti{q}+\Gamma\sum_{l=0}^{i-2}\omega^{i},
\end{flalign*}
so $\eqref{s1}_{2}$ yields that $\alpha_{i}\to \infty$. In these terms, \eqref{lip10} can be rearranged as
\begin{flalign}\label{lip12}
A_{i+1}\le \left[\frac{c2^{4i}(1+\kk_{i})^{2}}{(\tau_{2}-\tau_{1})^{2}}\right]^{\frac{2m\omega}{\alpha_{i+1}}}A_{i}^{\frac{\omega\alpha_{i}}{\alpha_{i+1}}},
\end{flalign}
for $c\equiv c(\texttt{data})$. Iterating \eqref{lip12} we obtain
\begin{flalign}\label{lip13}
A_{i+1}\le \left(\frac{c}{\tau_{2}-\tau_{1}}\right)^{\frac{4m}{\alpha_{i+1}}\sum_{l=1}^{i}\omega^{l}}\prod_{l=0}^{i-1}\left[2^{4(i-l)}(1+\kk_{i-l})^{2}\right]^{\frac{2m\omega^{l}}{\alpha_{i+1}}}A_{1}^{\frac{\omega^{i}\alpha_{1}}{\alpha_{i+1}}}.
\end{flalign}
Let us study the asymptotics of the various constants appearing in \eqref{lip13}. We have:
\begin{flalign*}
\lim_{i\to \infty}\frac{4m}{\alpha_{i+1}}\sum_{l=1}^{i}\omega^{l}=\frac{4m\omega}{\Gamma}, \qquad \lim_{i\to \infty}\frac{\omega^{i}\alpha_{1}}{\alpha_{i+1}}=\frac{m\ti{q}(\omega-1)}{\Gamma}
\end{flalign*}
and
\begin{flalign*}
\lim_{i\to \infty}\prod_{l=0}^{i-1}&\left[2^{8(i-l)}(1+\kk_{i-l})^{2}\right]^{\frac{2m\omega^{l}}{\alpha_{i+1}}}\nonumber \\
\le& \exp\left\{\frac{4m(\omega-1)}{\Gamma}\log\left(4\max\left\{2,\frac{\Gamma}{m(\omega-1)}\right\}\right)\left[\sum_{l=1}^{\left[\frac{1}{\log\omega}\right]+1}\omega^{-l}l+\frac{1+e^{-1}}{\log\omega}\right]\right\},
\end{flalign*}
where we also used that
\begin{flalign*}
\sum_{l=1}^{\infty}\omega^{-l}l\le \sum_{l=1}^{\left[\frac{1}{\log\omega}\right]+1}\omega^{-l}l+\frac{1+e^{-1}}{\log\omega}.
\end{flalign*} 
As
\begin{flalign}\label{lip14}
&\left(\mint_{Q_{\rr_{i+1}}}H_{j}(Du_{j})^{\alpha_{i+1}}\right)^{\frac{1}{\alpha_{i+1}}}\le A_{i+1}\nonumber \\
&\qquad \le\left(\frac{c}{\tau_{2}-\tau_{1}}\right)^{\frac{4m}{\alpha_{i+1}}\sum_{l=1}^{i}\omega^{l}}\prod_{l=0}^{i-1}\left[2^{4(i-l)}(1+\kk_{i-l})^{2}\right]^{\frac{2m\omega^{l}}{\alpha_{i+1}}}A_{1}^{\frac{\omega^{i}\alpha_{1}}{\alpha_{i+1}}},
\end{flalign}
we can pass to the limit in \eqref{lip14} for concluding that
\begin{flalign}\label{lip16}
&\nr{H_{j}(Du_{j})}_{L^{\infty}(Q_{\tau_{1}})}\le \frac{c}{(\tau_{2}-\tau_{1})^{\theta'}}\left(\mint_{Q_{\tau_{2}}}\left[1+H_{j}(Du_{j})^{m\ti{q}}\right] \ \dy\right)^{\frac{(\omega-1)}{\Gamma}}\nonumber \\
&\quad\le \frac{c}{(\tau_{2}-\tau_{1})^{\theta}}\left[1+\nr{H_{j}(Du_{j})}_{L^{\infty}(Q_{\tau_{2}})}^{\left(m\ti{q}-\frac{p}{2}\right)\frac{(\omega-1)}{\Gamma}}\right]\left(\mint_{Q_{\rr}}\left[1+H_{j}(Du_{j})^{\frac{p}{2}}\right] \ \dy\right)^{\frac{(\omega-1)}{\Gamma}},
\end{flalign}
with $c\equiv c(\texttt{data})$, $\theta'\equiv \theta'(n,p,q,d)$ and $\theta:=\theta'+(n+2)(\omega-1)\Gamma^{-1}$. Recalling the definition given in \eqref{r} and the restriction imposed in \eqref{pq}, it is easy to see that
\begin{flalign}\label{lip15}
\Gamma^{-1}\left(m\ti{q}-\frac{p}{2}\right)(\omega-1)<1.
\end{flalign}
In fact, verifying \eqref{lip15} is equivalent to check the validity of the following inequality
\begin{flalign*}
\ti{q}<\frac{p}{2m}+\frac{2}{\omega\ti{n}m},
\end{flalign*}
which is satisfied by means of \eqref{pq} and \eqref{tn}. So we can apply Young inequality with conjugate exponents $(b_{1},b_{2}):=\left(\frac{2\Gamma}{(2m\ti{q}-p)(\omega-1)},\frac{2\Gamma}{2\Gamma-(2m\ti{q}-p)(\omega-1)}\right)$ in \eqref{lip16} to end up with
\begin{flalign}\label{lip16.1}
&\nr{H_{j}(Du_{j})}_{L^{\infty}(Q_{\tau_{1}})}\le \frac{1}{2}\nr{H_{j}(Du_{j})}_{L^{\infty}(Q_{\tau_{2}})}\nonumber \\
&\quad+\frac{c}{(\tau_{2}-\tau_{1})^{\theta}}\left(\mint_{Q_{\rr}}\left[1+H_{j}(Du_{j})^{\frac{p}{2}}\right] \ \dy\right)^{\frac{(\omega-1)}{\Gamma}} +\frac{c}{(\tau_{2}-\tau_{1})^{\theta b_{2}}}\left(\mint_{Q_{\rr}}\left[1+H_{j}(Du_{j})^{\frac{p}{2}}\right] \ \dy\right)^{\frac{(\omega-1)b_{2}}{\Gamma}}\nonumber \\
&\quad \le \frac{1}{2}\nr{H_{j}(Du_{j})}_{L^{\infty}(Q_{\tau_{2}})}+\frac{c}{(\tau_{2}-\tau_{1})^{\theta b_{2}}}\left[1+\left(\mint_{Q_{\rr}}H_{j}(Du_{j})^{\frac{p}{2}} \ \dy\right)^{\frac{(\omega-1)b_{2}}{\Gamma}}\right],
\end{flalign} 
with $c\equiv c(\texttt{data})$. Now we apply Lemma \ref{l0} to \eqref{lip16.1} to conclude that
\begin{flalign}\label{lip17}
\nr{H_{j}(Du_{j})}_{L^{\infty}(Q_{\rr/2})}\le \frac{c}{\rr^{\beta_{1}}}\left[1+\left(\mint_{Q_{\rr}}H_{j}(Du_{j})^{\frac{p}{2}} \ \dy\right)^{\beta_{2}}\right],
\end{flalign}
for $c\equiv c(\texttt{data})$, $\beta_{1}:=\theta b_{2}$ and $\beta_{2}:=\frac{(\omega-1)b_{2}}{\Gamma}$.
\subsection{Proof of Theorem \ref{t1}} Let $\{u_{j}\}$ be the sequence built in Section \ref{inf}. As for each $j\in \N$, $u_{j}$ solves problem \eqref{pdjj}, which is driven by the nonlinear tensor $a_{j}(\cdot)$ defined in \eqref{lip3}, thus satisfying in particular \eqref{refregj}, and has boundary datum $f$ described by \eqref{ggg}, we deduce that the uniform energy bound \eqref{unibd} holds true. Hence, combining \eqref{unibd} with \eqref{lip17} we obtain that
\begin{flalign*}
\nr{H_{j}(Du_{j})}_{L^{\infty}(Q_{\rr/2})}\le\frac{c}{\rr^{\beta}}\left[\nr{Df}_{L^{r}(\Omega_{T})}^{r}+\nr{\partial_{t}f}^{p'}_{L^{p'}(0,T;W^{-1,p'}(\Omega))}+1\right],
\end{flalign*}
with $\beta:=\beta_{1}+(n+2)\beta_{2}$ and $c\equiv c(\texttt{data})$. Whenever $(t_{1},t_{2})\Subset (0,T)$ and $\tilde{\Omega}\Subset \Omega$ is open, a standard covering argument and the content of the above display render that
\begin{flalign}\label{lip18}
\nr{Du_{j}}_{L^{\infty}(\tilde{\Omega}\times (t_{1},t_{2}))}\le c(\texttt{data},\mathcal{C}_{f},\dist(\tilde{\Omega},\partial \Omega),t_{1},T-t_{2}).
\end{flalign}
Estimates \eqref{unibd} and \eqref{lip18} in turn imply that there exists a function $u\in L^{p}(0,T,W^{1,p}(\Omega))$ with gradient $Du\in L^{\infty}_{\loc}(0,T;L^{\infty}_{\loc}(\Omega,\mathbb{R}^{n}))$ so that
\begin{flalign}\label{lip19}
\begin{cases}
\ u_{j}\rightharpoonup u \quad &\mbox{in} \ \ L^{p}(0,T;W^{1,p}(\Omega))\\
\ Du_{j}\rightharpoonup^{*} Du\quad &\mbox{in} \ \ L^{\infty}_{\loc}(0,T;L^{\infty}_{\loc}(\Omega,\mathbb{R}^{n}))\\
\ u_{j}=f\quad &\mbox{on} \ \ \partial_{par}\Omega.
\end{cases}
\end{flalign}
In particular, by \eqref{lip18}, $\eqref{lip19}_{2}$ and weak$^{*}$-lower semicontinuity we have 
\begin{flalign}\label{lip27}
\nr{Du}_{L^{\infty}(\tilde{\Omega}\times (t_{1},t_{2}))}\le c(\texttt{data},\mathcal{C}_{f},\dist(\tilde{\Omega},\partial \Omega),t_{1},T-t_{2}).
\end{flalign}
Such information is not sufficient to pass to the limit as $j\to \infty$ in \eqref{wfj}, therefore we shall prove that $u_{j}$ admits some fractional derivative in space and in time which is controllable uniformly with respect to $j\in \N$. Concerning the fractional derivative in space, we can use \emph{verbatim} the same argument leading to \eqref{fr7}-\eqref{fr9} to deduce that
\begin{flalign*}
Du_{j}\in L^{p}_{\loc}(0,T;W^{\varsigma,p}_{\loc}(\Omega,\mathbb{R}^{n}))\quad \mbox{for all} \ \ \varsigma\in \left(0,\min\left\{1,\frac{2}{p}\right\}\right)
\end{flalign*}
with
\begin{flalign}\label{lip24}
\nr{u_{j}}_{L^{p}(t_{1},t_{2};W^{1+\varsigma}(\tilde{\Omega}))}\le c(\texttt{data},\varsigma,\mathcal{C}_{f},t_{1},T-t_{2},\dist(\tilde{\Omega},\partial \Omega)).
\end{flalign}
On the other hand, we cannot borrow the corresponding estimates for the fractional derivative in time of the $u_{j}$'s developed in \emph{Step 6} of Section \ref{high}: the constant appearing on the right-hand side of \eqref{fr9.1} depends on $\ti{\mu}^{-1}$ and, since now $\ti{\mu}\equiv \mu_{j}$, it may blow up in the limit as $j\to\infty$ if $\mu=0$. Therefore we shall follow a different path, see \cite[Section 9]{dumi1} for the case $q=p=2$. Let $0<t_{1}<\hat{t}_{1}<\hat{t}_{2}<t_{2}<T$ and $\tilde{h}>0$ be so that $0<\tilde{h}<\frac{\min\{\hat{t}_{1}-t_{1},t_{2}-\hat{t}_{2},1\}}{1000}$. Using the forward Steklov average to reformulate \eqref{wfj} we obtain, for a.e. $t\in (t_{1},t_{2})$,
\begin{flalign}\label{lip20}
\int_{\Omega}\left[\partial_{t}[u_{j}]_{\ti{h}}\varphi+[a_{j}(x,t,Du_{j})]_{\ti{h}}\cdot D\varphi\right] \ \dy=0\quad \mbox{for all} \ \ \varphi\in C^{\infty}_{c}(\Omega).
\end{flalign}
Since $\partial_{t}[u_{j}]_{\ti{h}}=\ti{h}^{-1}\ti{\tau}_{\ti{h}}u_{j}$, we can rearrange \eqref{lip20} as
\begin{flalign*}
\int_{\Omega}\left[\frac{\ti{\tau}_{\ti{h}}u_{j}}{\ti{h}}\varphi+[a_{j}(x,t,Du_{j})]_{\ti{h}}\cdot D\varphi\right] \ \dy=0.
\end{flalign*}
Modulo regularization, by \eqref{sss}, in the above display we can pick $\varphi:=\eta^{2}\ti{\tau}_{\ti{h}}u_{j}$ with 
\begin{flalign*}
\eta\in C^{\infty}_{c}(\tilde{\Omega})\quad \mbox{so that}\quad \nr{D\eta}_{L^{\infty}(\tilde{\Omega}}\le \frac{4}{\dist(\tilde{\Omega},\partial \Omega)},
\end{flalign*}
and integrate over the interval $(\hat{t}_{1},\hat{t}_{2})$ to get
\begin{flalign}\label{lip21}
h^{-1}\int_{\hat{t}_{1}}^{\hat{t}_{2}}\int_{\Omega}&\snr{\ti{\tau}_{\ti{h}}u_{j}}^{2}\eta^{2} \ \dx\ds=-\int_{\hat{t}_{1}}^{\hat{t}_{2}}\int_{\Omega}[a_{j}(x,t,Du_{j})]_{\ti{h}}\cdot\left[\eta^{2}\ti{\tau}_{\ti{h}}Du_{j}+2\ti{\tau}_{\ti{h}}u_{j}\eta D\eta \right] \ \dx\ds.
\end{flalign}
Recall that, for any function $w\in L^{1}(\ti{\Omega}\times (t_{1},t_{2}))$ there holds that
\begin{flalign*}
\int_{\hat{t}_{1}}^{\hat{t}_{2}}\int_{\ti{\Omega}}\snr{w_{\ti{h}}} \ \dx\ds\le \int_{\hat{t}_{1}-\ti{h}}^{\hat{t}_{2}+\ti{h}}\int_{\tilde{\Omega}}\snr{w} \ \dx\ds\le \int_{t_{1}}^{t_{2}}\int_{\tilde{\Omega}}\snr{w} \ \dx\ds,
\end{flalign*}
therefore, by $\eqref{refregj}_{1}$, \eqref{lip18}, H\"older and Young inequalities we estimate
\begin{flalign}\label{lip22}
&\left| \ \int_{\hat{t}_{1}}^{\hat{t}_{2}}\int_{\Omega}[a_{j}(x,t,Du_{j})]_{\ti{h}}\cdot\left[\eta^{2}\ti{\tau}_{\ti{h}}Du_{j}+2\ti{\tau}_{\ti{h}}u_{j}\eta D\eta \right] \ \dx\ds \ \right|\nonumber \\
&\quad \le 2\left(\sup_{\tilde{\Omega}\times (t_{1},t_{2})}\snr{Du_{j}}\right)\left(\int_{t_{1}}^{t_{2}}\int_{\tilde{\Omega}}a_{j}(x,t,Du_{j}) \ \dx\ds\right)\nonumber \\
&\quad+\frac{\ti{h}^{-1}}{2} \int_{\hat{t}_{1}}^{\hat{t}_{2}}\int_{\tilde{\Omega}}\eta^{2}\snr{\ti{\tau}_{\ti{h}}u_{j}}^{2} \ \dx \ds+\ti{h}\nr{D\eta}_{L^{\infty}(\tilde{\Omega}}^{2}\int_{t_{1}}^{t_{2}}\int_{\tilde{\Omega}}\snr{a_{j}(x,t,Du_{j})}^{2} \ \dx \dx\nonumber \\
&\quad \le \frac{\ti{h}^{-1}}{2} \int_{\hat{t}_{1}}^{\hat{t}_{2}}\int_{\tilde{\Omega}}\eta^{2}\snr{\ti{\tau}_{\ti{h}}u_{j}}^{2} \ \dx+c,
\end{flalign}
with $c\equiv c(\texttt{data},\mathcal{C}_{f},\dist{\tilde{\Omega},\partial \Omega},t_{1},T-t_{2})$. Merging \eqref{lip21} and \eqref{lip22} we end up with
\begin{flalign*}
\limsup_{\ti{h}\to \infty}\left(\ti{h}^{-1}\int_{\hat{t}_{1}}^{\hat{t}_{2}}\int_{\tilde{\Omega}}\snr{\ti{\tau}_{\ti{h}}u_{j}}^{2}\ \dx\ds\right)\le c(\texttt{data},\mathcal{C}_{f},\dist{\tilde{\Omega},\partial \Omega},t_{1},T-t_{2}),
\end{flalign*}
which, being $\hat{t}_{1},t_{1},\hat{t}_{2},t_{2}$ arbitrary, and since we can repeat exactly the same procedure for the backward Steklov average of $u_{j}$, we get
\begin{flalign*}
u_{j}\in W^{\iota,2}_{\loc}(0,T;L^{2}_{\loc}(\Omega))\quad \mbox{for all} \ \ \iota\in \left(0,\frac{1}{2}\right)
\end{flalign*}
and 
\begin{flalign}\label{lip23}
\nr{u_{j}}_{W^{\iota,2}(0,T;L^{2}(\tilde{\Omega}))}\le c(\texttt{data},\iota,\mathcal{C}_{f},\dist{\tilde{\Omega},\partial \Omega},t_{1},T-t_{2}).
\end{flalign}
From \eqref{lip24} and \eqref{lip23} we deduce that
\begin{flalign*}
\{u_{j}\} \ \mbox{is bounded uniformly w.r.t. $j\in \N$ in} \ W^{\iota,2}_{\loc}(0,T;L^{2}_{\loc}(\Omega))\cap L^{p}_{\loc}(0,T;W^{1+\varsigma,p}_{\loc}(\Omega))
\end{flalign*}
for all $\iota\in \left(0,\frac{1}{2}\right)$, $\varsigma\in \left(0,\min\left\{1,\frac{2}{p}\right\}\right)$, thus we can apply Lemma \ref{al} with $a_{1}= p$, $a_{2}=2$, $\sigma=\iota$, $X=W^{1+\varsigma,p}_{\loc}(\Omega)$, $B=W^{1,\min\{2,p\}}_{\loc}(\Omega)$, $Y=L^{2}_{\loc}(\Omega)$, to obtain a (non-relabelled) subsequence $\{u_{j}\}$ so that
\begin{flalign}\label{lip25}
u_{j}\to u\quad \mbox{in} \ \ L^{\min\{p,2\}}_{\loc}(0,T;W^{1,\min\{p,2\}}_{\loc}(\Omega)).
\end{flalign}
Combining $\eqref{lip19}_{2}$, \eqref{lip25} and \eqref{lip27} we get
\begin{flalign}\label{lip26}
Du_{j}\to Du\quad \mbox{in} \ \ L^{s}_{\loc}(0,T;L^{s}_{\loc}(\Omega,\mathbb{R}^{n}))\quad \mbox{for all} \ \ s\in (1,\infty),
\end{flalign}
therefore we can pass to the limit in \eqref{wfj} to deduce that $u$ satisfies \eqref{cv6}. Moreover, repeating \emph{Step 8} of Section \ref{high} we finally see that Definition \ref{d.1} is satisfied, therefore $u$ is a solution of problem \eqref{pdd} and, recalling also \eqref{lip27} we obtain $\eqref{t1.1}_{1}$. Once $\eqref{t1.1}_{1}$ is available, we can repeat the same procedure leading to \eqref{lip23} (with $a(\cdot),u$ replacing $a_{j}(\cdot),u_{j}$) to obtain $\eqref{t1.4}$. Furthermore, by \eqref{lip26}, \eqref{lip27} and \eqref{unibd} we can pass to the limit for $j\to \infty$ in \eqref{lip0} with $g\equiv 1$ and, after a standard covering argument, get $\eqref{t1.1}_{2}$. Finally, combining $\eqref{lip19}_{2}$ and \eqref{lip26} with \eqref{lip17} we obtain \eqref{t1.5}. The proof of Theorem \ref{t1} is complete.

\end{document}